%% file: main.tex
\documentclass[twocolumn]{autart}    

\usepackage{graphicx, caption, subcaption}
\usepackage{cite}
\usepackage{amsmath,amssymb,amsfonts}
\usepackage{textcomp}
\usepackage{xcolor}
\usepackage{algorithm}
\usepackage{siunitx,amsmath}
\usepackage{algpseudocode}
\usepackage{amsthm}
\newtheorem{assumption}{Assumption}
\newtheorem{lemma}{Lemma}
\newtheorem{definition}{Definition}
\newtheorem{remark}{Remark}
\newtheorem{theorem}{Theorem}
\newtheorem{proposition}{Proposition}
\newtheorem{corollary}{Corollary}

\newcommand{\stirlingii}{\genfrac{\langle}{\rangle}{0pt}{}}

\usepackage{easyReview}
\usepackage{fancyhdr}
\usepackage{enumitem}

\newcounter{reviewer}
\newcounter{point}[reviewer]
\usepackage{bibentry}
\renewcommand{\thepoint}{Q\,\thereviewer.\arabic{point}} 

\begin{document}

\begin{frontmatter}

\title{Gradient-tracking Based Differentially Private Distributed Optimization with Enhanced Optimization Accuracy\thanksref{footnoteinfo}} 

\thanks[footnoteinfo]{The work was supported in part by the National Science Foundation under
Grants ECCS-1912702, CCF-2106293, CCF-2215088, and CNS-2219487. This paper was not presented at any IFAC 
meeting. Corresponding author Yongqiang Wang. Tel. +1 864-656-5923. }

\author[USA]{Yu Xuan}\ead{yxuan@clemson.edu},    
\author[USA]{Yongqiang Wang}\ead{yongqiw@clemson.edu}               

\address[USA]{Department of Electrical and Computer Engineering, Clemson University, Clemson, SC 29634, USA}  

\begin{keyword}                           
Distributed optimization; Gradient tracking; Differential privacy.               
\end{keyword}                             

\begin{abstract}                          
Privacy protection has become an increasingly pressing requirement in distributed optimization. However, equipping distributed optimization with differential privacy, the state-of-the-art privacy protection mechanism, will unavoidably compromise optimization accuracy. In this paper, we propose an algorithm to achieve rigorous $\epsilon$-differential privacy in gradient-tracking based distributed optimization with enhanced optimization accuracy. {More specifically, to suppress the influence of differential-privacy noise, we propose a new robust gradient-tracking based distributed optimization algorithm that allows both stepsize and the variance of injected noise to vary with time. Then, we establish a new analyzing approach that can characterize the convergence of the gradient-tracking based algorithm under both constant and time-varying stespsizes. To our knowledge, this is the first analyzing framework that can treat gradient-tracking based distributed optimization under both constant and time-varying stepsizes in a unified manner. More importantly, the new analyzing approach gives a much less conservative analytical bound on the stepsize compared with existing proof techniques for gradient-tracking based distributed optimization.} We also theoretically characterize the influence of differential-privacy design on the accuracy of distributed optimization, which reveals that inter-agent interaction has a significant impact on the final optimization accuracy. Numerical simulation results confirm the theoretical predictions. 
\end{abstract}

\end{frontmatter}

\input{1_intro}
\input{2_preliminery}
\input{3_0_algorithm}
\input{3_1_convgent}

\input{3_2_DP}
\input{4_opt}
\input{5_simulation}

\section{Conclusions}

{In this paper, we have proposed a new differential-privacy approach for gradient-tracking based distributed optimization with enhanced optimization accuracy. By incorporating flexibilities in stepsize and noise variances, our approach achieves much better optimization accuracy than existing differential-privacy solutions for distributed optimization under the same privacy budget. To analyze {the} convergence of the proposed algorithm, we have also established a new proof technique that can address both constant and time-varying stepsizes. To our knowledge, this is the first analyzing framework that can address both time-varying and constant stepsizes for gradient-tracking based distributed optimization \textbf{in a unified framework}}. More importantly, the new proof technique gives a much less conservative analytical bound on the stepsize than those given in existing convergence analyzing approaches. We have also characterized the influence of interaction coupling weight on optimization accuracy. Numerical simulation results confirm the theoretical predictions.






\bibliographystyle{unsrt}        
\bibliography{main}           



\appendix

\input
\input{61_appendix1}
\input{62_appendix2}
\input{64_appendix4}


\newpage
\end{document}

%% file: 1_intro.tex
\section{Introduction}

We consider the problem of distributed optimization, where multiple agents cooperatively optimize a global aggregated objective function that is the sum of these agents' individual objective functions. Moreover, we assume that every agent can only share information with its immediate neighbors and there does not exist any central server to coordinate the communication or computation. Mathematically, the problem can be formulated as follows: 
\vspace{-2mm}
\begin{equation} \label{eq:DistOpt1}
    \min_{x\in\mathbb{R}^r} f\left(x\right)=\frac{1}{n}\sum_{i=1}^n f_i\left(x\right),
\end{equation}
where $f_i:\mathbb{R}^r\mapsto\mathbb{R},i\in[n]\triangleq \{1,\,2,\,\cdots,n\}$ is the local objective function only known to the $i^{th}$ agent, and $x\in\mathbb{R}^r$ is the decision or optimization variable. 

The study of the above distributed optimization problem can be traced back to the 1980s\cite{tsitsiklis1984problems}. In recent years, interests in such a problem have been revived by numerous emerging applications, including, for example, cooperative control\cite{yang2019survey}, distributed sensing \cite{bazerque2009distributed}, multi-agents systems\cite{raffard2004distributed}, sensor networks \cite{zhang2017distributed}, and distributed learning\cite{tsianos2012consensus}. To date, multiple solutions have been proposed for distributed optimization, with typical examples including gradient-based methods \cite{pu2021distributed,nedic2009distributed,shi2015extra,xu2017convergence,qu2017harnessing,xin2018linear,wang2022gradient}, distributed alternating direction method of multipliers \cite{shi2014linear,zhang2018admm}, and distributed Newton methods\cite{wei2013distributed}. Among these approaches, gradient based methods stand out due to their remarkable efficiency in computational complexity and storage requirement.

However, most existing gradient based solutions for problem~(\ref{eq:DistOpt1}) require participating agents to share explicit decision variables, which is unacceptable when the decision variables carry sensitive information of participating agents. For example, in the multi-robot rendezvous problem, the decision variable contains the position of participating robots, which should be kept private in unfriendly environments \cite{zhang2018admm,ruan2019secure,wang2019privacy,ridgley2020private,altafini2020system}.  In sensor-network based information fusion, the decision variables of participating agents should also be kept private as otherwise participating agents' location could be easily inferred \cite{zhang2018admm,huang2015differentially,burbano2019inferring,kia2015dynamic}. In distributed  machine learning applications, sharing decision variables or gradients has the risk of disclosing sensitive features of training data \cite{yan2012distributed,gade2018privacy}. In fact, recent studies in \cite{wang2022decentralized} show that without a privacy mechanism in place, an adversary can use shared information to precisely recover
the raw data used for training (pixel-wise accurate for images
and token-wise matching for texts).

To ensure the privacy of participating agents in distributed optimization, plenty of approaches have been reported. For example, partially homomorphic encryption has been employed in our work \cite{zhang2018admm,zhang2018enabling} as well as others' work \cite{freris2016distributed,lu2018privacy,hadjicostis2020privacy}  to enable privacy protection in distributed optimization. This approach, however, will incur heavy communication and computational overhead. {Time or spatially correlated ``structured'' noise based approaches have also been proposed by us \cite{wang2019privacy,gao2022algorithm,gao2022dynamics} as well as others \cite{manitara2013privacy,mo2016privacy,he2018preserving,altafini2020system,yan2012distributed,lou2017privacy,gade2018private}   to enable privacy protection. Because the added noise can cancel out (either temporally or spatially), this approach can ensure the accuracy of distributed optimization. However,  since the noises injected are correlated, the enabled strength of privacy protection is limited. In fact, to protect the privacy of one agent, such approaches usually require the agent to have at least one neighbor that does not share information with the adversary, which may be hard to satisfy for general multi-agent systems.} Due to its wide applicability and simplicity, differential privacy (DP) has become the de facto standard for privacy protection. Recently, differential privacy has also been applied to distributed optimization using independent additive noises \cite{huang2015differentially,cortes2016differential,xiong2020privacy,ding2021differentially,nozari2016differentially,han2022differentially,chen2021differential,wang2022tailoring}. However, these independent noises also unavoidably reduce the accuracy of distributed optimization. 

Leveraging a recently proposed  technique to track the cumulative gradient \cite{pu2020robust}, we incorporate time-varying stepsize and noise variance in gradient-tracking based distributed optimization and propose a distributed optimization framework that can reduce the influence of differential-privacy noise on optimization accuracy while achieving rigorous $\epsilon$-differential privacy. We also propose a new analyzing approach that can characterize the convergence under both time-varying and constant stepsizes for gradient-tracking based distributed optimization. To the best of our knowledge, this is the first analyzing approach able to achieve this goal in a unified framework. {It is worth noting that \cite{pu2021distributed} considers gradient-tracking based distributed optimization under both constant and time-varying stepsizes. However, it uses different proof techniques for the two cases. In addition, our considered decaying stepsizes are more general than the decaying stepsize considered in \cite{pu2021distributed}, which is a special case of ours.} By  characterizing the influence of differential-privacy noises on the accuracy of distributed optimization, we discover that inter-agent interaction has a significant impact on the final optimization error.

The main contributions are summarized as follows: 
     {1) We propose a gradient-tracking based distributed optimization approach that has guaranteed $\epsilon$-differential privacy but reduced influence of differential-privacy design on optimization accuracy. Compared with existing differential-privacy results for distributed optimization (both the result for static-consensus based distributed optimization in \cite{huang2015differentially} and the result for gradient-tracking based distributed optimization in \cite{ding2021differentially}), simulation results show that the proposed approach performs much better in optimization accuracy under the same privacy level (privacy budget); 
     2) To characterize the convergence of the proposed algorithm, we establish a new analyzing framework that can address both constant and time-varying stepsizes. To the best of our knowledge, this is the first proof technique that is able to characterize the convergence of gradient-tracking based distributed optimization under both constant and time-varying stepsizes in a unified framework. More importantly, the new analyzing approach gives a much less conservative analytical bound on stepsize compared with existing proof techniques for gradient-tracking based distributed optimization. For instance, for the example used {in} the numerical simulation, our new convergence proof gives a stepsize bound that is 100 times larger than the one given in \cite{pu2020robust}: our analysis yields an upper bound on stepsize $\alpha$ as $1.2\times 10^{-3}$, while the stepsize bound from the analysis in \cite{pu2020robust} is $1.1\times 10^{-5}$;
    3) We prove that under decreasing stepsizes, our algorithm ensures $\epsilon$-differential privacy in the infinite time horizon. To our knowledge, our algorithm is the first to achieve this goal for gradient tracking based distributed optimization under general objective functions. (Note that \cite{ding2021differentially} requires adjacent objective functions to have identical Lipschitz constants and convexity parameters); }
    4) We theoretically characterize the influence of the inter-agent coupling weights on the final optimization accuracy in the presence of differential-privacy noises.

The remainder of the paper is organized as follows: Sec. 2 presents some preliminaries and the problem formulation. Sec. 3 proposes a gradient-tracking based distributed optimization algorithm and rigorously characterizes the convergence performance and privacy strength. Sec. 4  systematically investigates the influence of inter-agent coupling weights on the final optimization error. Sec. 5 uses numerical simulation results { to show the advantages of the proposed approach over existing counterparts.}

%% file: 2_preliminery.tex
\section{Preliminaries and problem formulation} 
\subsection{Notations}
\vspace{-0.3cm} 
We use a lower case letter in bold or a capital letter, such as $\boldsymbol{x}$ or $X$, to denote a matrix, and a lower case letter, such as $x$, to denote a scalar or a column vector. 
The inner product of two vectors is denoted as $\langle x,y\rangle$, and the inner product of two matrices for $X,Y\in\mathbb{R}^{n\times r}$ is defined as $\langle X,Y\rangle\triangleq \sum_{i=1}^n\langle x_i,y_i\rangle$, where $x_i$ and $y_i$ represent the $i^{th}$ row of $X$ and $Y$, respectively.
The 2-norm for a vector $x\in\mathbb{R}^n$ is defined as $\|x\|_2=\sqrt{\sum_{i=1}^n {x_i^2}}$, where $x_i$ is the $i^{th}$ element. Furthermore, given an arbitrary vector norm $\|\cdot\|$ on $\mathbb{R}^n$, we define the matrix norm for any $\boldsymbol{x}\in\mathbb{R}^{n\times r}$ as $\|\boldsymbol{x}\|\triangleq\left\|\ \left[\|x^{(1)}\|,\|x^{(2)}\|,\dots,\|x^{(r)}\|\right]\ \right\|_2$, where $x^{(i)}$ is the $i^{th}$ column of $\boldsymbol{x}$ and $\|\cdot\|_2$ represents the 2-norm. 
$\rho(\cdot)$ denotes the spectral radius of a  matrix. We use $[n]$ to represent the set $\{1, 2, \dots, n\}$. {Finally, we define $d_I^2\triangleq\|I-\frac{\boldsymbol{1}\boldsymbol{1}^T}{n}\|_2^2$ and  $\rho_w \triangleq \rho\left(W-\frac{\boldsymbol{1}\boldsymbol{1}^T}{n}\right)$.}

We use a graph to describe the inter-agent interaction. More specifically, an undirected graph is denoted as a pair $\mathcal{G} = (\mathcal{V},\mathcal{E})$, where $\mathcal{V}$ represents the set of agents, and $\mathcal{E}\subseteq\mathcal{V}\times\mathcal{V}$ represents the set of undirected edges. $(j,i) \in\mathcal{E}$ represents the edge connecting the $i^{th}$ and $j^{th}$ agents. We describe the inter-agent interaction using an undirected graph induced by a nonnegative symmetric matrix $W=[w_{ij}]\in\mathbb{R}^{n\times n}$, i.e., $\mathcal{G}_W=([n],\mathcal{E}_W)$, and $(j,i)\in\mathcal{E}_W$ if and only if $w_{ij}>0$. $\mathcal{N}_i$ denotes the set of neighboring agents of $i$, i.e.,  $\mathcal{N}_i=\{j\neq i| (i,j)\in\mathcal{G}_W\}$.

\subsection{Distributed optimization} \label{sec:pre}
\vspace{-0.3cm}

For the convenience of differential-privacy analysis, we describe the distributed optimization problem $\mathcal{P}$ in (1) by three parameters ($\mathcal{X},f,\mathcal{G}_W$):

\vspace{-2mm}
(i) {$\mathcal{X}\subseteq\mathbb{R}^r$} is the domain of optimization;
\\
{
(ii) $f(x)\triangleq\frac{1}{n}\sum_{i=1}^n f_i(x)$ with $f_i:\mathbb{R}^r\mapsto\mathbb{R}$;}
\\
(iii) $\mathcal{G}_W$ specifies the undirected communication graph induced by a weight matrix $W$.

We rewrite problem~(\ref{eq:DistOpt1}) as an equivalent problem below:
\vspace{-2mm}
\begin{equation} \label{eq:DistOpt2}
    \begin{split}
        \min_{x_1,x_2,\dots,x_n \in \mathbb{R}^r} & \frac{1}{n}\sum_{i=1}^{n}f_i\left(x_i\right)\\
        \text{subject to } & x_1=x_2=\dots=x_n,
    \end{split}
    \vspace{-2mm}
\end{equation}
where $x_i$ is agent $i$'s local copy of the decision variable.

{We make the following assumptions on objective functions in (\ref{eq:DistOpt2}):}
\vspace{0.2cm}
\begin{assumption} \label{ass:muL}
    Each $f_i$ is $\mu$-strongly convex with an $L$-Lipschitz continuous gradient, i.e., for any {$x,x'\in\mathbb{R}^r$},
    \vspace{-2mm}
    \begin{equation}
    \begin{split}
        \langle \nabla{f_i\left(x\right)}-\nabla{f_i\left(x'\right)},x-x'\rangle&\geq \mu\|x-x'\|^2,\\
        \|\nabla{f_i\left(x\right)}-\nabla{f_i\left(x'\right)}\|&\leq L\|x-x'\|.
    \end{split}
    \end{equation}
\end{assumption}
\vspace{-7mm}
{Under Assumption~\ref{ass:muL}, Problem~(\ref{eq:DistOpt1}) has a unique solution $x^\ast\in\mathbb{R}^r$.}
\vspace{0.2cm}
\begin{assumption} \label{ass:C}
    The gradients of all local objective functions are bounded on $\mathcal{X}$. Without loss of generality, we denote $C>0$ as the bound such that for any $x\in\mathcal{X}$ and any $f_i\in\mathcal{S}$, $\|\nabla f_i(x)\|_2\leq C$.
\end{assumption}
\vspace{-2mm}
Note that the assumption of bounded gradients can be relaxed in practice by using gradient clipping.

%% file: 3_0_algorithm.tex
\section{A gradient-tracking based differentially private distributed optimization algorithm }
For a network of $n$ agents with a coupling matrix $W$, the proposed algorithm is summarized in Algorithm 1:
\begin{algorithm} 
\caption{Gradient-tracking based differentially private distributed optimization} \label{alg:DP}
\begin{algorithmic}
    \State \hspace{-4mm} {\textbf{Parameters: } stepsizes $\alpha$ and  $\gamma_k$; noise factor $\beta_k$; inter-agent coupling weights  $w_{ij}>0$ for $\forall (i,j)\in\mathcal{E}_W$; Laplace DP noises $\eta_{i,k}\triangleq [\eta_{i1,k}, \eta_{i2,k}, \dots, \eta_{ir,k}]^T$ with $\eta_{ij,k}\sim Lap(b_\eta)$ and $\xi_{i,k}\triangleq [\xi_{i1,k}, \xi_{i2,k}, \dots, \xi_{ir,k}]^T$ with $\xi_{ij,k}\sim Lap(b_\xi)$;} 
    \State \hspace{-4mm} {\textbf{Initialization: }Each agent $i\in[n]$ initializes with arbitrary states $x_{i,0}$, $s_{i,0}$ in $\mathcal{X}\subseteq\mathbb{R}^r$};
    \State \hspace{-4mm} \textbf{for } $k=0, 1, \dots,$ \textbf{do} 
    \State a) Agent $i$ injects DP noises {$\beta_k\eta_{i,k}$ and $\beta_k\xi_{i,k}$} on $s_{i,k}$ and $x_{i,k}$, respectively, and sends $w_{ji}(s_{i,k}+{\beta_k\eta_{i,k}})$, and $w_{ji}(x_{i,k}+{\beta_k\xi_{i,k}})$  to each agent $j\in\mathcal{N}_i$;
    \State b) After receiving $w_{ij}(x_{j,k}+{\beta_k\xi_{j,k}})$ and $w_{ij}(s_{j,k}+{\beta_k\eta_{j,k}})$ from each agent $j\in\mathcal{N}_i$, each agent $i$ updates its states as follows:
    {
    \vspace{-2mm}
    \begin{equation*}
    \begin{split}
        \hspace{-2cm}{s_{i,k+1}} \hspace{-0.8mm}=& {w_{ii}}s_{i,k}\hspace{-.8mm}+\hspace{-3.5mm}\sum_{j=1,j\neq i}^n\hspace{-3mm}{w_{ij}}(s_{j,k} \hspace{-.5mm}+\hspace{-.5mm} \beta_k\eta_{j,k}) \hspace{-.5mm}+\hspace{-.5mm} \gamma_k\nabla{f_i(x_{i,k})},\\
        {x_{i,k+1}} \hspace{-0.8mm}=& {w_{ii}}x_{i,k}\hspace{-.8mm}+\hspace{-3.5mm}\sum_{j=1,j\neq i}^n\hspace{-3mm}{w_{ij}}(x_{j,k} \hspace{-.5mm}+\hspace{-.5mm} \beta_k\xi_{j,k}) \hspace{-.5mm}-\hspace{-.5mm} \alpha(s_{i,k+1} \hspace{-.5mm}-\hspace{-.5mm} s_{i,k}).
    \end{split}
    \end{equation*}}
    \State \hspace{-4mm} \textbf{end for}
\end{algorithmic}
\end{algorithm}

Note that inspired by \cite{pu2020robust}, in our Algorithm 1, we let individual agents track a cumulative gradient $s_i$ instead of the gradient in commonly used gradient-tracking approaches. This can avoid noise from accumulating in the estimation of the global gradient. {But it is worth noting that different from \cite{pu2020robust}{,} which only considers constant stepsize and noise variances, we allow the stepsize and noise variances to be time-varying, which is crucial to {achieving} differential privacy in the infinite time horizon, as detailed later in Sec. 3.2.}

From Algorithm 1, it can be seen that to achieve $\epsilon$-differential privacy, each agent $i$ injects noises $\eta_{i,k}\in\mathbb{R}^{r}$ and $\xi_{i,k}\in\mathbb{R}^r$ when it exchanges information with its neighbors.  
Denote $\boldsymbol{\eta}_k\triangleq[\eta_{1,k},\eta_{2,k},\dots,\eta_{n,k}]^T\in \mathbb{R}^{n\times r}$ and $\boldsymbol{\xi}_k\triangleq[\xi_{1,k},\xi_{2,k},\dots,\xi_{n,k}]^T\in \mathbb{R}^{n\times r}$. We make the following assumption on the differential-privacy noise:

\vspace{0.2cm}
\begin{assumption} \label{ass:noise}
    The elements of random sequences $\{\boldsymbol{\eta}_k\}$ and $\{\boldsymbol{\xi}_k\}$ are independent along iteration $k$, across agents $i$, and state dimension $j$. $\boldsymbol{\eta}_k$ and $\boldsymbol{\xi}_k$ have zero means and bounded variances, i.e., $\mathbb{E}[\boldsymbol{\eta_k}]=\mathbb{E}[\boldsymbol{\xi_k}]=\boldsymbol{0}$,  $\mathbb{E}[\|\boldsymbol{\eta}_k\|^2]\leq \sigma^2_\eta$, and  $\mathbb{E}[\|\boldsymbol{\xi}_k\|^2]\leq \sigma^2_\xi$ for some  $\sigma_\eta$, $\sigma_\xi>0$. 
\end{assumption}

Defining the augmented $x_{i}$, $s_{i}$, and $\nabla f_i(x_i)$ as $\boldsymbol{x}\triangleq[x_1,x_2,\dots,x_n]^T\in \mathbb{R}^{n\times r}$,  $\boldsymbol{s}\triangleq[s_1,s_2,\dots,s_n]^T\in \mathbb{R}^{n\times r}$, and $\nabla \hspace{-.5mm}F(\boldsymbol{x}) \hspace{-.5mm}\triangleq\hspace{-.5mm} [\nabla \hspace{-.5mm}f_1(x_1),\hspace{-.5mm}\nabla \hspace{-.5mm}f_2(x_2),\dots,\hspace{-.5mm}\nabla \hspace{-.5mm}f_n(x_n)]^T\in \mathbb{R}^{n\times r} $, we can rewrite Algorithm 1 in the following form:
\vspace{-2mm}
\begin{equation} \label{eq:algdecay}
    {
    \begin{split}
    \boldsymbol{s}_{k+1} & = W{\boldsymbol{s}_k} + {W_o}\beta_k\boldsymbol{\eta}_k+\gamma_k\nabla{F\left(\boldsymbol{x}_k\right)},
     \\
    \boldsymbol{x}_{k+1} & = W\boldsymbol{x}_k + {W_o}\beta_k\boldsymbol{\xi}_k -\alpha\left(\boldsymbol{s}_{k+1}-\boldsymbol{s}_k\right),
    \end{split}}
\end{equation}
\hspace{-1mm}
where $W$ is the weight matrix, and ${W_o}$ is obtained by replacing all diagonal entries of $W$ with zeros while keeping all other entries unchanged. 
{We make the following assumption on the coupling matrix $W$:}
\vspace{0.2cm}
\begin{assumption} \label{ass:network}
    The graph $\mathcal{G}_W$ is undirected and connected, i.e., there is a path between any two agents. The matrix $W\in \mathbb{R}^{n\times n}$ is non-negative, doubly-stochastic, and symmetric, i.e., $w_{ij}\geq0$,  $W\boldsymbol{1}=\boldsymbol{1},$  $\boldsymbol{1}^TW=\boldsymbol{1}^T$, and $W=W^T$. In addition, $w_{ii}>0$ for all  $i\in[n].$
\end{assumption} 

To make connections with conventional gradient-tracking based algorithms, define $\boldsymbol{y}_k$ as $\boldsymbol{y}_k\triangleq \boldsymbol{s}_{k+1}-\boldsymbol{s}_k$. Then one can obtain
\vspace{-2mm}
\begin{equation} \label{eq:DPpp2}
    \begin{split}
        \boldsymbol{y}_{k+1} & = W {\boldsymbol{y}_{k}}+\Tilde{\nabla}F(\boldsymbol{x}_{k+1})-\Tilde{\nabla}F(\boldsymbol{x}_{k}),\\
        \boldsymbol{x}_{k+1} & = W\boldsymbol{x}_k + {W_o}{\beta_k}\boldsymbol{\xi}_k -\alpha\boldsymbol{y}_k,
    \end{split}
\end{equation}
\hspace{-1mm}
{where $\Tilde{\nabla}F(\boldsymbol{x}_k) = \gamma_k\nabla F(\boldsymbol{x}_k)+{W_o}\beta_k\boldsymbol{\eta}_k$}.
One can verify that the average $\boldsymbol{y}_k$, i.e, $\frac{\boldsymbol{1}^T\boldsymbol{y}_k}{n}$, tracks the global gradient: \hspace{-4mm}
\begin{equation} \label{eq:1Tyk}
   \boldsymbol{1}^T\boldsymbol{y}_k = \boldsymbol{1}^T\Tilde{\nabla}F\left(\boldsymbol{x}_k\right), \forall k. 
   \vspace{-3mm}
\end{equation}

%% file: 3_1_convgent.tex
\vspace{-2mm}
\subsection{Convergence analysis} \label{sec:conv}
\vspace{-1mm}
In this section, we consider time-varying stepsize $\gamma_k$ and noise factor $\beta_k$. More specifically, we set $\gamma_k=\frac{\gamma}{(m+k)^p}$ and $\beta_k=\frac{1}{(m+k)^q}$ with $p\geq0$, $q\geq0$, and $m>0$. For the convenience of analysis, we define  $\bar{x}_k\triangleq\frac{1}{n}\boldsymbol{1}^T{\boldsymbol{x}_k}$, $\bar{y}_k\triangleq\frac{1}{n}\boldsymbol{1}^T{\boldsymbol{y}_k}$, and
\vspace{-3mm}
\begin{equation*}
    \left\{\begin{array}{rcl}
        h_k & \triangleq\frac{1}{n}\boldsymbol{1}^T\nabla F\left(\boldsymbol{x}_k\right) & = \frac{1}{n} \sum^n_{i=1}\nabla f_i\left(x_{i,k}\right) \\
        g_k & \triangleq\frac{1}{n}\boldsymbol{1}^T\nabla F\left(\boldsymbol{1}\bar{x}_k\right) & = \frac{1}{n} \sum^n_{i=1}\nabla f_i\left(\bar{x}_k\right) 
    \end{array}.\right.
    \vspace{-3mm}
\end{equation*} 
{Denoting $v^{T}\triangleq\boldsymbol{1}^{T}W_o$, one can obtain the dynamics of $\bar{x}_{k+1}$, $\boldsymbol{x}_{k+1}-\boldsymbol{1}\bar{x}_{k+1}$, and $\boldsymbol{y}_{k+1}-\boldsymbol{1}\bar{y}_{k+1}$
as:
\vspace{-2mm}
\begin{subequations}
    \vspace{-2mm}
    \begin{equation} \label{eq:barx}
    \bar{x}_{k+1} 
        \hspace{-.5mm} = \hspace{-.5mm} \bar{x}_k \hspace{-.5mm}-\hspace{-.5mm} \alpha\gamma_k g_k \hspace{-.5mm}-\hspace{-.5mm} \alpha\gamma_k(h_k \hspace{-.5mm} -\hspace{-.5mm} g_k) \hspace{-.5mm}-\hspace{-.5mm} \frac{\alpha}{n}v^T\hspace{-.5mm}\beta_k\boldsymbol{\eta}_k
        \hspace{-.5mm}+\hspace{-.5mm} \frac{1}{n}v^T\hspace{-.5mm}\beta_k\boldsymbol{\xi}_k,
    \end{equation}
    \vspace{-5mm}
    \begin{equation} \label{eq:dynx}
    \begin{array}{ll}
    \hspace{-8mm}\boldsymbol{x}_{k+1} \hspace{-1mm}-\hspace{-.5mm} \boldsymbol{1}\bar{x}_{k+1} \hspace{-1mm}
        & \hspace{-.8mm }= \hspace{-.8mm}\left(\hspace{-1mm}W \hspace{-.8mm} - \hspace{-.8mm} \frac{\boldsymbol{1}\boldsymbol{1}^{\hspace{-.5mm}T}}{n}\hspace{-.5mm}\right)\hspace{-1mm}(\boldsymbol{x}_k\hspace{-.8mm}-\hspace{-.8mm}\boldsymbol{1}\bar{x}_k) \hspace{-.5mm}-\hspace{-.5mm} \alpha\hspace{-1mm} \left(\hspace{-1mm}I \hspace{-.8mm}-\hspace{-.5mm} \frac{\boldsymbol{1}\boldsymbol{1}^{\hspace{-.5mm}T}}{n}\hspace{-.5mm}\right)\hspace{-.5mm}\\
        & \hspace{3mm}\times\hspace{-.5mm}(\boldsymbol{y}_k \hspace{-.5mm}-\hspace{-.5mm} \boldsymbol{1}\bar{y}_k)
        \hspace{-.5mm}+\hspace{-.5mm} \left(I - \frac{\boldsymbol{1}\boldsymbol{1}^T}{n}\right)\hspace{-.5mm}{W_o}\boldsymbol{\xi}_k,
    \end{array}
    \end{equation}
    \vspace{-2mm}
    \begin{equation} \label{eq:dyny}
    \begin{array}{ll}
    \hspace{-3mm}\boldsymbol{y}_{k+1} \hspace{-1mm}-\hspace{-.5mm}\boldsymbol{1}\bar{y}_{k+1}\hspace{-1.5mm}
        &= \hspace{-1mm}\left(W \hspace{-.5mm}-\hspace{-.5mm} \frac{\boldsymbol{1}\boldsymbol{1}^T}{n}\right)\hspace{-.5mm}\left(\boldsymbol{y}_k \hspace{-.5mm}-\hspace{-.5mm} \boldsymbol{1}\bar{y}_k\right)\\
        & \hspace{3.5mm} + \hspace{-.5mm}
        \left(\hspace{-.5mm}I \hspace{-.5mm}-\hspace{-.5mm} \frac{\boldsymbol{1}\boldsymbol{1}^T}{n}\hspace{-.5mm}\right)\hspace{-1mm}\left(\hspace{-.5mm}\Tilde{\nabla}\hspace{-.5mm}F(\boldsymbol{x}_{k+1}) \hspace{-.5mm}-\hspace{-.5mm} \Tilde{\nabla}\hspace{-.5mm}F(\boldsymbol{x}_k)\hspace{-.5mm}\right)\hspace{-1mm}.
    \end{array}
    \end{equation}
\end{subequations}
}
{
\vspace{1mm}
\begin{lemma} \label{lemmaIneq}
Under Assumptions 1, 3, 4, and $p\geq 0$, when $\alpha\gamma_0<\frac{2}{\mu+L}$, we have
\vspace{-4mm}
\begin{equation} \label{eq:linear}
        \left[  \begin{matrix} \vspace{-1mm}
        \mathbb{E}[\|\bar{x}_{k+1}\hspace{-0.5mm}-x^\ast\|_2^2]\\
        \mathbb{E}[\|\boldsymbol{x}_{k+1}\hspace{-0.5mm}-\hspace{-0.5mm}\boldsymbol{1}\bar{x}_{k+1}\|_2^2]\\
        \mathbb{E}[\|\boldsymbol{y}_{k+1}\hspace{-0.5mm}-\hspace{-0.5mm}\boldsymbol{1}\bar{y}_{k+1}\|_2^2]
        \end{matrix}\right]\hspace{-1mm}\leq
        \hspace{-0.5mm}A_k\hspace{-0.5mm}\left[  \begin{matrix}
        \mathbb{E}[\|\bar{x}_{k}\hspace{-0.5mm}-\hspace{-0.5mm}x^\ast\|_2^2]\\
        \mathbb{E}[\|\boldsymbol{x}_{k}\hspace{-0.5mm}-\hspace{-0.5mm}\boldsymbol{1}\bar{x}_{k}\|_2^2]\\
        \mathbb{E}[\|\boldsymbol{y}_{k}\hspace{-0.5mm}-\hspace{-0.5mm}\boldsymbol{1}\bar{y}_{k}\|_2^2]
        \end{matrix}\right]\hspace{-1.5mm} +\hspace{-.5mm} B_k,
        \vspace{-2mm}
    \end{equation}
    where the inequality is taken in a component-wise manner, and $A_k$ and $B_k$ are given by
    \vspace{-3mm}
    \begin{equation} \label{eq:A_k}
    \begin{split}
        & A_k \hspace{-.8mm} =\hspace{-1.5mm}
        \left[ \hspace{-.5mm}\begin{matrix}
        1\hspace{-.8mm}-\hspace{-.5mm}\Tilde{\alpha}_k\mu  & 0 & {A_k^{31}} \\
        \Tilde{\alpha}_k a_{12} & \frac{1+\rho_w^2}{2} & {A_k^{32}} \\
         0 & \alpha^2a_{23} & \frac{1+\rho_w^2}{2} \hspace{-.8mm}+\hspace{-.8mm}  \Tilde{\alpha}_{k+1}^2a_{33}
        \end{matrix}\right]^{\hspace{-.5mm}T} 
        \hspace{-3mm}, 
        B_k  \hspace{-.8mm}=\hspace{-1.5mm} \left[ \begin{matrix}
        \beta_k^2 b_{1} \\
        \beta_k^2 b_{2} \\
        {B_k^3}
        \end{matrix}\right]{\hspace{-1mm},}
        \vspace{-4mm}
    \end{split}
    \vspace{-5mm}
    \end{equation}
    with  
    $\tilde{\alpha}_k\triangleq\alpha\gamma_k$, $A_k^{31} \hspace{-.5mm}= \gamma_{k+1}^2\Tilde{\alpha}_k^2a_{31} + (\gamma_k \hspace{-.3mm}-\hspace{-.3mm} \gamma_{k+1})^2a_{34}$, 
    $A_k^{32} \hspace{-.5mm}= \gamma_{k+1}^2a_{32}$
    $ + (\gamma_k \hspace{-.3mm}-\hspace{-.3mm} \gamma_{k+1})^2a_{35}$, $B_k^3 \hspace{-.5mm} = \beta_k^2 b_{31} \hspace{-.5mm}+\hspace{-.5mm} (\gamma_k\hspace{-.5mm}-\hspace{-.5mm}\gamma_{k+1})^2 b_{32}$, 
    $a_{12} = \frac{L^2}{\mu n}(1+\Tilde{\alpha}_k\mu)$,    
    $a_{23} = \frac{(1+\rho_w^2)}{1-\rho_w^2}$,
    $a_{31} = \frac{32n L^4d_I^2}{1-\rho_w^2}$, \\
    $a_{32} = \frac{32d_I^2(L^2+\Tilde{\alpha}_k^2L^4)}{1-\rho_w^2}$,
    $a_{33} = \frac{16L^2d_I^2}{1-\rho_w^2}$,
    $a_{34}=\frac{16nL^2d_I^2}{1-\rho_w^2}$,
    $a_{35}=\frac{16L^2d_I^2}{1-\rho_w^2}$,
    $b_1 = \frac{\alpha^2\|v\|_2^2}{n^2}\sigma_\eta^2 + \frac{\|v\|^2_2}{n^2}\sigma_\xi^2$, 
    $b_2 = d_I^2\|{W_o}\|_2^2\sigma_\xi^2$,
    $b_{31} = \hspace{-1mm}\frac{4d_I^2\sigma_\eta^2}
    {1-\rho_w^2}\hspace{-1mm}\left(\hspace{-.5mm}\frac{4\Tilde{\alpha}_{k+1}^2L^2\|v\|_2^2}{n}
    \hspace{-.5mm}+\hspace{-.5mm}  (1 \hspace{-.5mm}+\hspace{-.5mm} \Tilde{\alpha}_{k+1} L )\|{W_o}\|_2^2\right)$\\
    $+ \frac{4\gamma_{k+1}^2L^2d_I^2}{1-\rho_w^2} \|{W_o}\|_2^2\sigma_\xi^2$, and 
    $b_{32} = \frac{8d_I^2\sum_{i=1}^n\|f_i(x^\star)\|_2^2}{1-\rho_w^2}$.
\end{lemma}}
\vspace{-5mm}
\begin{proof}
See Appendix A.
\end{proof}
\vspace{-4mm}
{Based on Lemma~\ref{lemmaIneq}, we can obtain the following convergence result for Algorithm~\ref{alg:DP}: }
\vspace{2mm}
{
\begin{theorem} \label{thmConv}
    Under Assumptions 1, 3, and 4, there exists some $m>0$ such that under $\gamma_k=\gamma/(m+k)^p$ and  $\beta_k= 1/(m+k)^q$, the errors of Algorithm~\ref{alg:DP} satisfy
    \begin{enumerate}[itemsep=2pt,topsep=0pt,parsep=0pt, leftmargin=14pt, itemindent=0pt] \vspace{-2mm}
 \item[1)]  ${\left\{\begin{array}{cc}   
        \mathbb{E}   [\|\bar{x}_{k}-x^\ast\|_2^2]                  \leq\frac{\mathcal{O}(1)}{(m+k)^{2q}}\\
        \mathbb{E}[\|\boldsymbol{x}_{k} -\boldsymbol{1}\bar{x}_{k}\|_2^2]
            \leq\frac{\mathcal{O}(1)}{(m+k)^{2q}}\\
         \mathbb{E}[\|\boldsymbol{y}_{k}-\boldsymbol{1}\bar{y}_{k}\|_2^2]
            \leq\frac{\mathcal{O}(1)}{(m+k)^{2q}}
        \end{array}\right.,}$
        \\
        \\
        with $\mathcal{O}(1)$ denoting some constant independent of $k$, if $p=0$ (i.e., $\gamma_k=\gamma$), $q>0$, and $\alpha$ and $\gamma$ satisfy $\alpha\gamma<\min\left\{\frac{2}{\mu+L}, \frac{1-\rho_w^2}{4\sqrt{2}Ld_I},\sqrt{\frac{2}{c_{42}+\sqrt{c_{42}^2+4c_{41}}}}\right\}$, with $c_{41}=\frac{64(1+\rho_w^2)(3\mu+L)L^6d_I^2}{(1-\rho_w^2)^4(\mu+L)}$ and $c_{42}=\frac{64(1+\rho_w^2)(\mu+L+2L^2)L^2d_I^2}{(1-\rho_w^2)^4(\mu+L)}$;
\item[2)]  $\left\{\begin{array}{cc} 
        \mathbb{E}[\|\bar{x}_{k}-x^\ast\|_2^2]
                      \leq \frac{\mathcal{O}(1)}{(m+k)^{\min\{2q-p,2p\}}}\\
        \mathbb{E}[\|\boldsymbol{x}_{k} -\boldsymbol{1}\bar{x}_{k}\|_2^2]
                    \leq\frac{\mathcal{O}(1)}{(m+k)^{2\min\{q,p\}}}\\
        \mathbb{E}[\|\boldsymbol{y}_{k}-\boldsymbol{1}\bar{y}_{k}\|_2^2]                     \leq\frac{\mathcal{O}(1)}{(m+k)^{2\min\{q,p\}}}
        \end{array}\right.,$
        \\
        \\
        with $\mathcal{O}(1)$ denoting some constant independent of $k$, if $0< p\leq 1$ and $q> p$, and $\alpha$ and $\gamma$ satisfy $\alpha\gamma>\frac{\min\{2q-p,2p\}}{\mu}$;
\item[3)]  ${\left\{\begin{array}{cc}                     \mathbb{E}[\|\bar{x}_{k}\hspace{-0.5mm}-            x^\ast\|_2^2]\leq \mathcal{O}(1) +                         \frac{\mathcal{O}(1)}{(m+k)^p}\\
       \mathbb{E}[\|\boldsymbol{x}_{k} -\boldsymbol{1}\bar{x}_{k}\|_2^2]
                    \leq\frac{\mathcal{O}(1)}{(m+k)^{2\min\{q,p\}}}\\
        \mathbb{E}[\|\boldsymbol{y}_{k}-\boldsymbol{1}\bar{y}_{k}\|_2^2]
                    \leq\frac{\mathcal{O}(1)}{(m+k)^{2\min\{q,p\}}}
        \end{array}\right.,}$
        \\
        \\
        with $\mathcal{O}(1)$ denoting some constant independent of $k$, if $p > 1$ and $q\geq\frac{p}{2}$.
    \end{enumerate}
\end{theorem}
\vspace{-4mm}
\begin{proof}
    See Appendix B.
\end{proof}
}
\begin{corollary} \label{remark:decayrate}
    When $p=0$ and $q\rightarrow\infty$, for $0\leq k \leq q$, the error evolutions of Algorithm~\ref{alg:DP} satisfy \vspace{-3mm}
\begin{equation}
\vspace{-3mm}
    \begin{array}{l} \label{eq:geodecay}
        \mathbb{E}[\|\bar{x}_{k}\hspace{-0.5mm}-x^\ast\|_2^2]
        \leq e^{-\frac{4k}{2+\ln (1/a)} }\mathcal{O}(1),\\
        \mathbb{E}[\|\boldsymbol{x}_{k} -\boldsymbol{1}\bar{x}_{k}\|_2^2]
        \leq e^{-\frac{4k}{2+\ln (1/a)} } \mathcal{O}(1),\\
        \mathbb{E}[\|\boldsymbol{y}_{k}-\boldsymbol{1}\bar{y}_{k}\|_2^2]
        \leq e^{-\frac{4k}{2+\ln (1/a)} } \mathcal{O}(1),
    \end{array}
\end{equation}
with some constants $a\in(0,1)$ and $\mathcal{O}(1)$ independent of $k$
if
$\alpha$ and $\gamma$ satisfy \\$\alpha\gamma\hspace{-.5mm}<\hspace{-.5mm}\min\hspace{-.5mm}\left\{\hspace{-.5mm}\frac{2}{\mu+L}, \frac{1-\rho_w^2}{4\sqrt{2}Ld_I}, \hspace{-,5mm}\sqrt{\frac{2}{c_{42}+\sqrt{c_{42}^2+4c_{41}}}}\hspace{-.5mm}\right\}$, with $c_{41}=\frac{64(1+\rho_w^2)(3\mu+L)L^6d_I^2}{(1-\rho_w^2)^4(\mu+L)}$
and $c_{42}=\frac{64(1+\rho_w^2)(\mu^2+2\mu L+5L^2)L^2d_I^2}{(1-\rho_w^2)^4(\mu+L)^2}$.
\end{corollary}
\begin{proof}
    {
    From the proof of Theorem~\ref{thmConv}, the condition of $m$ for case 1) in Theorem~\ref{thmConv} is given as
    \begin{equation} \label{eq:cond0p}
    \begin{array}{l}
    \left(\frac{m}{m+1}\right)^{2q} > \max\left\{ 1-\alpha\gamma\mu, \frac{1+\rho_w^2}{2}+{\alpha^2\gamma^2a_{33}}\right\}, \\
    1 \hspace{-.5mm}>\hspace{-.5mm} {\frac{\alpha^5\gamma^5a_{23}a_{31}L^2(3\mu \hspace{-.1mm}+\hspace{-.2mm} L)}{\mu n(\mu\hspace{-.1mm}+\hspace{-.2mm}L)c_1c_2c_3}} \hspace{-.5mm}+\hspace{-.5mm} {\frac{\alpha^2\hspace{-.2mm}\gamma^2 32a_{23}d_I^2\hspace{-.2mm}L^{\hspace{-.2mm}2}\hspace{-.2mm}({\mu^2\hspace{-.3mm}+\hspace{-.2mm}2\mu L \hspace{-.2mm}+\hspace{-.2mm} 5L^2})}{(1\hspace{-.2mm}-\hspace{-.2mm}\rho_w^2)(\mu\hspace{-.1mm}+\hspace{-.2mm}L)^2c_2c_3}},
    \end{array}
    \vspace{-2mm}
    \end{equation}
    with $c_1 = \left(\frac{m}{m+1}\right)^{2q} +\alpha\gamma\mu - 1$, $c_2 = 2\left(\frac{m}{m+1}\right)^{2q} - (1+\rho_w^2)$, $c_3 = \left(\frac{m}{m+1}\right)^{2q}-\frac{1+\rho_w^2}{2}-{\alpha^2\gamma^2a_{33}}$, and constants $a_{23}$, $a_{31}$, and $a_{33}$ given in Lemma~\ref{lemmaIneq}. 
Define $a\triangleq\left(\frac{m}{m+1}\right)^{2q}$ with $m>0$. Reversely, we have $m=\frac{a^{1/2q}}{1-a^{1/2q}}$, $0<a<1$. Rewrite the constraints of $m$ in (\ref{eq:cond0p}) in terms of $a$ as
\vspace{-2mm}
    \begin{equation*} 
    \begin{array}{ll}
    & a >
    \max\left\{ 1-\alpha\gamma\mu, \frac{1+\rho_w^2}{2}+{\alpha^2\gamma^2a_{33}}\right\}, \\
    & (a- \frac{1+\rho_w^2}{2})(a - \frac{1+\rho_w^2}{2}-{\alpha^2\gamma^2a_{33}}) \\
    & \hspace{3mm}>\hspace{-1.3mm} {\frac{\alpha^5\gamma^5a_{23}a_{31}L^2(3\mu + L)}{\mu n(\mu + L)(a  + \alpha\gamma\mu- 1)}} \hspace{-.5mm}+\hspace{-.5mm} {\frac{\alpha^2\gamma^232a_{23}d_I^2L^2({\mu^2+2\mu L + 5L^2})}{(1-\rho_w^2)(\mu+L)^2}}.
    \end{array}
    \end{equation*}
The existence of $m$ ensures the existence of $a$.
Next, we show that $m$ is in the order of $\mathcal{O}(cq)$ as $q$ goes to infinity, with $c$ denoting some positive constant. 
Denoting $z\triangleq \frac{1}{2q}$, by L'Hopital's Rule, we have
\vspace{-2mm}
\begin{equation*}
\begin{array}{ll}
    \lim_{q\rightarrow\infty} \frac{a^{1/2q}}{cq(1-a^{1/2q})} & = \lim_{z\rightarrow0} \frac{\frac{2}{c}za^{z}}{(1-a^{z})}\\
    & = \lim_{z\rightarrow0} \frac{2}{c} \frac{a^z + z a^z\ln a}{-a^z\ln a }=\frac{2}{c}\left(-\frac{1}{\ln a}\right).
\end{array}
\end{equation*}
When $c=\frac{-\ln a}{2}$, we have $\lim_{q\rightarrow\infty} \frac{a^{1/2q}}{cq(1-a^{1/2q})}=1$, i.e., $\lim_{q\rightarrow\infty}m=\frac{-\ln a}{2}q$.}

{
Note that {when $p=0$}, the error convergence bound decays at the rate $\frac{1}{(m+k)^{2q}}$, and the ratio of the decaying rate between two adjacent iterations $k$ and $k+1$ satisfies
\vspace{-2mm}
\begin{equation*}
    \frac{\frac{1}{(m+k)^{2q}}}{\frac{1}{(m+k+1)^{2q}}} = \left(1+\frac{1}{m+k}\right)^{2q}.
\end{equation*}
For $0\leq k \leq q$, as $q$ tends to infinity we have
\vspace{-2mm}
\begin{equation*}  
\begin{split}
    \lim_{q\rightarrow\infty}\left(1+\frac{1}{m+k}\right)^{2q}>\lim_{q\rightarrow\infty}\left(1+\frac{1}{m+q}\right)^{2q}
    = e^{\frac{4}{2-\ln a}}.
\end{split}
\end{equation*}
Namely, as $q\rightarrow\infty$, for all $0\leq k \leq q$, we have
\vspace{-2mm}
\begin{equation*}
    \frac{1}{(m+k+1)^{2q}}< e^{-\frac{4}{2+\ln (1/a)} } \frac{1}{(m+k)^{2q}}.
\end{equation*}
Therefore, when $p=0$, as $q\rightarrow\infty$, the error for Algorithm~\ref{alg:DP} reaches a geometric decaying rate $\mathcal{O}\left(e^{-\frac{4k}{2+\ln (1/a)} }\right)$ in  (\ref{eq:geodecay}).}
\end{proof}
\vspace{-4mm}
{
When the stepsize is set as a constant value, as in existing gradient-tracking based distributed optimization algorithms like \cite{pu2020robust}, Corollary~\ref{remark:decayrate} gives a much less conservative analytical bound on the stepsize:
\vspace{2mm}
\begin{corollary} \label{remark:alpha}
    Under $\gamma=1$ and $n >1$, the stepsize bound given in Corollary~\ref{remark:decayrate} is larger than the stepsize bound given in \cite{pu2020robust}.   
\end{corollary}}
\vspace{-4mm}
\begin{proof}
    {
Under coupling matrix $W$, the stepsize condition in \cite{pu2020robust} is\\
$\alpha<\min\left\{\frac{1}{\mu+L}, \frac{1-\rho_w^2}{4\sqrt{3}Ld_I}, \sqrt{\frac{2d_3}{d_2+\sqrt{d_2^2+4d_1d_3}}}\right\}$,
with 
$d_1 = \frac{48d_I^4L^6}{\mu(1-\rho_w^2)^2}$,
$d_2 = \frac{24d_I^4 L^2({2L^2} + \mu^2 n)(\|W-I\|_2^2 + 2)}{\mu(1-\rho_w^2)^2} + \frac{10L^4d_I^2}{\mu}$, and
$d_3\triangleq \frac{\mu n (1-\rho_w^2)^2}{18}$.}

{
Recall that for $p=0$, $q>0$, and $\gamma=1$, $\alpha$ satisfies
\\
$\alpha<\min\left\{\frac{2}{\mu+L}, \frac{1-\rho_w^2}{4\sqrt{2}Ld_I},\sqrt{\frac{2}{c_{42}+\sqrt{c_{42}^2+4c_{41}}}}\right\}
$, with $c_{41}$ and $c_{42}$ given in Corollary~\ref{remark:decayrate}.
To prove the corollary, we only need to prove
\vspace{-2mm}
\begin{equation}
    \sqrt{\frac{2}{c_{42}+\sqrt{c_{42}^2+4c_{41}}}}>\sqrt{\frac{2}{\frac{d_2}{d_3}+\sqrt{\frac{d_2^2}{d_3^2}+4\frac{d_1}{d_3}}}},
\end{equation}
with
$\frac{d_2}{d_3} = \frac{432d_I^4 L^2({2L^2} + \mu^2 n)(\|W-I\|_2^2 + 2)}{\mu^2n(1-\rho_w^2)^4} + \frac{180L^4d_I^2}{\mu^2 n (1-\rho_w^2)^2}$ and 
$\frac{d_1}{d_3} = \frac{864 d_I^4L^6}{\mu^2 n (1-\rho_w^2)^4}$. 
It is sufficient to show for $n>1$ and $\mu\leq 2$ that the following inequalities hold:
\begin{equation} \label{eq:paraineq1}
    c_{42}<\frac{d_2}{d_3},
\end{equation}
\vspace{-3mm}
\begin{equation} \label{eq:paraineq2}
    c_{42}^2+4c_{41}<\frac{d_2^2}{d_3^2}+4\frac{d_1}{d_3}.
\end{equation}
Given $1+\rho_w^2<2$, we have $c_{42}<\left(\frac{1}{2}+\frac{2L^2}{(\mu+L)^2}\right)\frac{256^2d_I^2}{(1-\rho_w^2)^4}$. Given $\frac{d_I^2}{n} = \frac{n-1}{n} \geq \frac{1}{2}$, we have $\frac{d_2}{d_3} > \left(1+\frac{2L^2}{\mu^2}\right)\frac{432d_I^2 L^2}{(1-\rho_w^2)^4}$. Since $L>0$, we always have (\ref{eq:paraineq1}). Similarly, given $1+\rho_w^2<2$, $\frac{d_I^2}{n} = \frac{n-1}{n} \geq \frac{1}{2}$, and $0<1-\rho_w^2<1$, we can prove (\ref{eq:paraineq2}). 
}
\end{proof}
\begin{remark} \label{remark:simalpha}
    Under the ridge regression problem with randomly generated parameters considered in numerical simulations in Sec. 5.2, our Corollary~\ref{remark:decayrate} yields an analytical upper bound on the stepsize $\alpha$ as $1.2 \times 10^{-3}$. For the same set of parameters, the convergence analysis in \cite{pu2020robust} gives an upper bound on the stepsize as $1.1 \times 10^{-5}$. Hence, our analysis can provide a much less conservative stepsize bound.
\end{remark}

%% file: 3_2_DP.tex
\subsection{Differential-privacy analysis} \label{sec:dp} 

We denote the $j^{th}$ element of $s_{i,k}$  ( $\eta_{i,k}$,  $\xi_{i,k}$, and $x_{i,k}$, respectively) as $s_{ij,k}$ (  $\eta_{ij,k}$, $\xi_{ij,k}$, and $x_{ij,k}$, respectively). In Algorithm 1, noises are injected in both the decision variables $\boldsymbol{x}_k$ and auxiliary variables $\boldsymbol{s}_k$. We denote the observed $\boldsymbol{x}_k$ and $\boldsymbol{s}_k$ as $\boldsymbol{m}_k\triangleq[m_{1,k},m_{2,k},\dots,m_{n,k}]^T\in\mathbb{R}^{n\times r}$ and $\boldsymbol{n}_k\triangleq[n_{1,k},n_{2,k},\dots,n_{n,k}]^T\in \mathbb{R}^{n\times r}$, respectively, where $m_{ij,k}\triangleq s_{ij,k}+{\beta_k}\eta_{ij,k}$ and $n_{ij,k}\triangleq x_{ij,k}+{\beta_k}\xi_{ij,k}$.

{To achieve differential privacy, Laplace noise\cite{dwork2006calibrating} and Gaussian noise \cite{dwork2014algorithmic,canonne2020discrete} are two common choices. {Since} privacy analysis is easier
under Laplace noise, we consider Laplace noise in the paper. 
\vspace{2mm}
\begin{assumption} \label{ass:noiseLap}
    $\{\boldsymbol{\eta}_k\}$ and $\{\boldsymbol{\xi}_k\}$ follow Laplace distribution, i.e., 
    \vspace{-4mm}
    \begin{equation} \label{eq:lap}
    \eta_{ij,k}\sim Lap\left(b_{\eta}\right) \text{ and } \xi_{ij,k}\sim Lap\left(b_{\xi}\right),
    \vspace{-3mm}
\end{equation}
where $Lap(b)$ denotes the Laplace distribution with probability density function $f_L(x)\triangleq \frac{1}{2b}e^{-|x|/b}$.
\end{assumption}
}
Recall that $Lap(b)$ has zero mean and variance $2b^2$. Also, for any $x,y\in\mathbb{R}$, we have $\frac{{f_L}(x)}{{f_L}(y)}\leq e^{\frac{|y-x|}{b}}$.

Given a problem $\mathcal{P}$, we describe an iterative distributed optimization algorithm as a mapping $Alg(\mathcal{P},\boldsymbol{x}(0)):\boldsymbol{x}(0)\in\mathcal{X}\mapsto\mathcal{O}\in\mathbb{O}$, where $\boldsymbol{x}(0)$ is the initial state, $\mathcal{O}$ denotes the observation sequence, and $\mathbb{O}$ denotes the set of all possible observation sequences. To be specific, for Algorithm~\ref{alg:DP}, the observation sequence is $\mathcal{O} = \{\boldsymbol{m}_0,\boldsymbol{n}_0,\boldsymbol{m}_1,\boldsymbol{n}_1,\dots\}$, and the initial state is $\boldsymbol{x}(0)=\{\boldsymbol{s}_0,\boldsymbol{x}_0\}$.
Moreover, given any $\mathcal{O}\in\mathbb{O}$, we denote $\mathcal{O}_{k} \triangleq \{\boldsymbol{m}_0,\boldsymbol{n}_0,\dots,\boldsymbol{m}_{k-1},\boldsymbol{n}_{k-1}\}$ as the subsequence of $\mathcal{O}$ with $\mathcal{O}_{0} \triangleq \emptyset$. Given any distributed optimization problem $\mathcal{P}$, observations $\mathcal{O}_k$, and initial states $\boldsymbol{x}(0)$ under Algorithm \ref{alg:DP}, we denote $\mathcal{A}(\mathcal{P},\mathcal{O}_{k},\boldsymbol{x}(0))$ as the internal states at iteration $k$, i.e., $\{\boldsymbol{x}_k$, $\boldsymbol{s}_k\}$.

Inspired by \cite{huang2015differentially} and \cite{dwork2010differential}, we define adjacency, $\epsilon$-differential privacy, and sensitivity for Algorithm~\ref{alg:DP} as follows: 
\vspace{0.2cm}
{
\begin{definition} \label{def1}
Two distributed optimization problems $\mathcal{P}$ and $\mathcal{P}'$ are adjacent if the following conditions hold:
\\
(i) $\mathcal{X}=\mathcal{X}'$ and $\mathcal{G}_W=\mathcal{G}_W'$, i.e., the domain of optimization and the communication graphs are identical; 
\\
(ii) there exists an $i\in[n]$ such that $f_i\neq f_i'$ and for all $j\neq i,f_j=f_j'$.
\end{definition}
}

From Definition \ref{def1}, it can be seen that two distributed optimization problems are adjacent if and only if one agent changes its local objective function while all other parameters remain the same.

\vspace{0.2cm}
\begin{definition} \label{def:dp}
    ($\epsilon$-differential privacy) For a given $\epsilon>0$, an iterative distributed algorithm solving problem~(\ref{eq:DistOpt2}) is $\epsilon$-differentially private if for any two adjacent problems $\mathcal{P}$ and $\mathcal{P}'$, any set of observation sequences $\mathbb{Y}\subseteq\mathbb{O}$, and any initial state $\boldsymbol{x}(0)$, we always have
    \vspace{-2mm}
    \begin{equation}
    \mathbb{P}\left(Alg(\mathcal{P},\boldsymbol{x}(0) ) \in \mathbb{Y}\right)
    \leq e^{\epsilon} \mathbb{P}\left(Alg(\mathcal{P}',\boldsymbol{x}(0) ) \in \mathbb{Y}\right),
    \end{equation}
    where the probability is taken over the randomness of iteration processes.
\end{definition} 
\vspace{2mm}
\begin{remark}
It is worth noting that different from \cite{huang2015differentially}, which defines DP based on the probability distribution of internal states, we follow the standard DP framework in \cite{dwork2010differential} and use the probability distribution of observations in Definition \ref{def:dp}.
\end{remark}
\vspace{-2mm}
$\epsilon$-differential privacy makes sure that an adversary can hardly identify a local objective function among all possible ones, even when all the observation sequences are disclosed. A smaller $\epsilon$ implies a higher privacy level.
{
\vspace{2mm}
\begin{definition} \label{defses}
    (Sensitivity) At each iteration $k\in[K]\triangleq\{0,1,\dots,K\}$, given any initial  state $\boldsymbol{x}(0)=\{\boldsymbol{x}_0,\boldsymbol{s}_0\}$, for adjacent distributed optimization problems $\mathcal{P}$ and $\mathcal{P}'$, the sensitivity in $\boldsymbol{s}$ and $\boldsymbol{x}$ are respectively defined as
    \vspace{-2mm}
    \begin{equation*} 
    \begin{split}
        & \Delta s_k \triangleq \sup_{\mathcal{O}\in\mathbb{O}}
        \hspace{-1mm}\sup_{\scriptsize \begin{array}{l} \vspace{-1mm}
        \boldsymbol{s}_k\hspace{-.5mm}\in\hspace{-.5mm}\mathcal{A}(\mathcal{P},\mathcal{O}_{k},\boldsymbol{x}(0)) \\
        \boldsymbol{s}_k'\hspace{-.5mm}\in\hspace{-.5mm}\mathcal{A}(\mathcal{P}',\mathcal{O}_{k},\boldsymbol{x}(0))\end{array}}\hspace{-5mm}\|\boldsymbol{s}_k-\boldsymbol{s}_k'\|_1,
        \end{split}
    \end{equation*}
    \vspace{-6mm}
    \begin{equation*}
        \begin{split}
        & \Delta x_k \triangleq \sup_{\mathcal{O}\in\mathbb{O}}
        \hspace{-1mm}\sup_{\scriptsize\begin{array}{l}\vspace{-1mm}
        \boldsymbol{x}_k\hspace{-.5mm}\in\hspace{-.5mm}\mathcal{A}(\mathcal{P},\mathcal{O}_{k},\boldsymbol{x}(0)) \\
        \boldsymbol{x}_k'\hspace{-.5mm}\in\hspace{-.5mm}\mathcal{A}(\mathcal{P}',\mathcal{O}_{k},\boldsymbol{x}(0))\end{array}}\hspace{-5mm}\|\boldsymbol{x}_k-\boldsymbol{x}_k'\|_1.
    \end{split}
    \vspace{-5mm}
    \end{equation*}
    \vspace{-4mm}
\end{definition}
It is well known that the differential privacy of an algorithm depends on its sensitivity. Specifically, we give the following proposition:
\vspace{2mm}
\begin{proposition} \label{prop:DP}
    Under Assumptions~\ref{ass:noise} and \ref{ass:noiseLap}, Algorithm~\ref{alg:DP} is $\epsilon$-differentially private for a given privacy budget $\epsilon>0$ if the following condition holds:
    \vspace{-2mm}
\begin{equation} \label{eq:epsilonCond0}
    \begin{split}
        \sum_{k\in[K]}\left(\frac{\Delta s_k}{\beta_k b_{\eta}}+\frac{\Delta x_k}{\beta_k b_{\xi}}\right)\leq {\epsilon}.
    \end{split}
    \vspace{-2mm}
\end{equation}
\vspace{-4mm}
\end{proposition}
\vspace{-5mm}
\begin{proof}
    In our algorithmic setting, the probability distribution of an observation sequence is the same as the probability distribution of the corresponding internal state sequence in \cite{huang2015differentially}. Therefore, following a derivation similar to Lemma 2 in [21], we can obtain the proposition.
\end{proof}
\vspace{-4mm}
The following theorem gives the privacy budget for $Algorithm$~\ref{alg:DP} under general stepsize $\gamma_k$ and noise factor $\beta_k$:
\vspace{3mm}
\begin{theorem} \label{thmDP}
When Assumptions 2-5 hold, given any finite number of iterations $K$, Algorithm 1 is $\epsilon$-differentially private for a given $\epsilon>0$ if the noise parameters $b_\eta$ and $b_\xi$ in Assumption~\ref{ass:noiseLap} satisfy
\vspace{-2mm}
\begin{equation} \label{eq:thmDP}
\begin{split}       {2\sqrt{r}C}\sum_{k=1}^{K}\sum_{t=0}^{k-1}\left(\frac{w_{ii}^{k-1-t}}{\beta_k b_\eta} + \frac{\alpha|c_{k,t}|}{\beta_k b_\xi} \right)\gamma_t
\leq \epsilon
\end{split}
\end{equation}
\vspace{-1mm}
for all $i\in[n]$,
where 
\vspace{-2mm}
\begin{equation} \label{eq:ckt0}
    c_{k,t}\triangleq w_{ii}^{k-2-t}((k-t-1)-(k-t)w_{ii}).
\end{equation}
\vspace{-1mm}
\end{theorem}
\vspace{-5mm}

\begin{proof}
Given $\boldsymbol{x}_0=\boldsymbol{x}'_0$,  $\boldsymbol{s}_0=\boldsymbol{s}'_0$, any observation sequence $\mathcal{O}$, and an arbitrary pair of adjacent problems $\mathcal{P}$ and $\mathcal{P}'$, denote $\Delta s_{i,k} \triangleq s_{i,k}-s'_{i,k}$, $\Delta x_{i,k} \triangleq x_{i,k}-x'_{i,k}$, and  $\Delta f_{i,k}\triangleq\nabla f_i(x_{i,k}) - \nabla f'_i(x'_{i,k})$ for $s_{i,k}, x_{i,k}\in\mathcal{A}(\mathcal{P},\mathcal{O}_k,\boldsymbol{x}(0))$ and $s'_{i,k}, x'_{i,k}\in\mathcal{A}(\mathcal{P}',\mathcal{O}_k,\boldsymbol{x}(0))$. For all $k\geq0$, we have
\vspace{-2mm}
\begin{equation*}
    \begin{split}
        & \Delta s_{i,k+1}
        =w_{ii}\Delta s_{i,k}+ \gamma_k\Delta f_{i,k},\\
        &\Delta x_{i,k+1}
        =w_{ii}\Delta x_{i,k} 
        + \alpha (1 - w_{ii})\Delta s_{i,k} - \alpha\gamma_k\Delta f_{i,k}.
    \end{split}
\end{equation*}
For all $i\in[n]$, we have $\Delta s_{i,0} = \Delta x_{i,0} = 0$. Based on the definition of adjacent problems $\mathcal{P}$ and $\mathcal{P}'$, for $i$ satisfying $f_i\neq f'_i$, we have the following relationship by induction for all $1\leq k \leq K$:
\vspace{-2mm}
\begin{equation} \label{eq:35}
\begin{split}
    \Delta s_{i,k} = \sum_{t=0}^{k-1}w_{ii}^{k-1-t}\gamma_t\Delta f_{i,t}, 
    \Delta x_{i,k} = \sum_{t=0}^{k-1} \alpha {c_{k,t}} \gamma_t\Delta f_{i,t
    },
    \end{split}
\end{equation}
where $c_{k,t}$ is defined in (\ref{eq:ckt0}).
For all $j$ satisfying $f_j=f'_j$, we always have $\Delta s_{j,k}=\Delta x_{j,k} = \boldsymbol{0}$.
\\
\\
Hence we have the following relationship:
\vspace{-2mm}
\begin{equation} \label{eq:39}
\begin{split}
    & \sum_{k=0}^K\hspace{-.5mm}\frac{\|\Delta s_{k}\|_1}{\beta_k}
    \hspace{-.5mm}=\hspace{-.5mm}\sum_{k=1}^K\hspace{-.5mm}\frac{\|\Delta s_{i,k}\|_1}{\beta_k}
    \hspace{-.5mm}\leq\hspace{-.5mm}
    \sum_{k=1}^{K}\hspace{-.5mm}\sum_{t=0}^{k-1}\hspace{-.5mm}\frac{w_{ii}^{k-1-t}}{\beta_k}\gamma_t\|\Delta f_{i,t}\|_1,
    \\
    &\sum_{k=0}^K\hspace{-.5mm}\frac{\|\Delta x_{k}\|_1}{\beta_k}
    \hspace{-.5mm}=\hspace{-.5mm}\sum_{k=0}^K\hspace{-.5mm}\frac{\|\Delta x_{i,k}\|_1}{\beta_k}
    \hspace{-.5mm}\leq\hspace{-.5mm} \sum_{k=1}^{K}\hspace{-.5mm}\sum_{t=0}^{k-1}\hspace{-.5mm}\alpha\frac{|c_{k,t}|}{\beta_k}\gamma_t\|\Delta f_{i,t}\|_1.
\end{split}
\end{equation}
Under Assumption \ref{ass:C}, we have
\vspace{-2mm}
\begin{equation} \label{eq:40}
    \|\Delta f_{i,k}\|_1 = \|\nabla f_i(x_{i,k})-\nabla f'_i(x'_{i,k})\|_1 \leq 2\sqrt{r}C.
\end{equation}
Using Proposition~\ref{prop:DP}, combining (\ref{eq:39}) and (\ref{eq:40}) leads to the result in (\ref{eq:thmDP}).
\end{proof}

\vspace{-2mm}
{When $\gamma_k$ and $\beta_k$ decrease with time, we can prove that Algorithm~\ref{alg:DP} can ensure $\epsilon$-differential privacy even when the number of iterations tends to infinity:}

\vspace{2mm}
\begin{corollary} \label{corDPinf}
    When Assumptions 2-5 hold, Algorithm~\ref{alg:DP} is $\epsilon$-differentially private even when the number of iterations tends to infinity if $\gamma_k=\frac{\gamma}{(m+k)^p}$ and $\beta_k=\frac{1}{(m+k)^q}$ with $m>0$, $p\geq0$, $q<p-2$, and $b_\eta$ and $b_\xi$ in Assumption~\ref{ass:noiseLap} satisfy 
    \vspace{-2mm}  
    \begin{equation} \label{eq:thmDPinf}
    \begin{split}
    \frac{2\sqrt{r}C\gamma {\rm P}_{{\lceil p \rceil}}\hspace{-.5mm}(w_{ii})}{m^p(1 \hspace{-.5mm}-\hspace{-.5mm}w_{ii})^{{\lceil p \rceil}+1}}\hspace{-1.5mm}
    \left( \hspace{-1mm}
    \frac{\sum_{k=m+1}^{\infty}\hspace{-.5mm}\frac{1}{k^{p-q}}}{b_\eta w_{ii}^{m}} 
    \hspace{-.8mm}
    + \hspace{-.8mm}
    \frac{\alpha\hspace{-.5mm}\sum_{k=m+1}^{\infty}\hspace{-.5mm}\frac{1}{k^{p-q-1}}}{b_\xi w_{ii}^{m+1}} \hspace{-1.5mm}\right)\hspace{-1.5mm}
    \leq \hspace{-.8mm} \epsilon
    \end{split}
    \vspace{-8mm}
    \end{equation}  
    for all $i\in[n]$, where ${\lceil p \rceil}$ is the least integer that is no less than $p$ and ${\rm P}_{{\lceil p \rceil}}(w_{ii})$ is the Eulerian polynomial of order ${\lceil p \rceil}$, i.e., 
    \vspace{-4mm}
    \begin{equation} \label{eq:pn}
        {\rm P}_{{\lceil p \rceil}}(w_{ii}) \triangleq \left\{ \begin{array}{cc}
          \sum_{j=0}^{{\lceil p \rceil}-1} \stirlingii{{\lceil p \rceil}}{j} w_{ii}^j   &  {\lceil p \rceil}\geq 1\\
    1         & {\lceil p \rceil}=0
        \end{array}\right.,
    \end{equation}
    with $\stirlingii{{\lceil p \rceil}}{j}$ being the entry $j$ in row ${\lceil p \rceil}$ of the triangle of Eulerian numbers \cite{edgar2018staircase} (note that $\stirlingii{{\lceil p \rceil}}{j}$ is always $1$ for $j=0$). (Note that under the given conditions for $p$, $q$, and $m$, $\sum_{k=m+1}^{\infty}\hspace{-.5mm}\frac{1}{k^{p-q}}$ and $\sum_{k=m+1}^{\infty}\hspace{-.5mm}\frac{1}{k^{p-q-1}}$ are always finite. The Eulerian polynomial ${\rm P}_{{\lceil p \rceil}}(w_{ii})$ is also finite under any finite $p$.)
\end{corollary}
\vspace{-5mm}
\begin{proof}
    In light of the condition $0<w_{ii}<1$ implied from Assumption \ref{ass:network}, the coefficient $c_{k,t}$ defined in (\ref{eq:ckt0}) satisfies
\vspace{-2mm}
\begin{equation} \label{eq:ckt}
\begin{split}
    |c_{k,t}| \hspace{-.5mm}=\hspace{-.5mm}   w_{ii}^{k \hspace{-.2mm}-\hspace{-.2mm} 2 \hspace{-.2mm}-\hspace{-.2mm} t}|(k \hspace{-.5mm}-\hspace{-.5mm} t \hspace{-.5mm}-\hspace{-.5mm} 1)(1 \hspace{-.5mm}-\hspace{-.5mm} w_{ii}) \hspace{-.5mm}-\hspace{-.5mm} w_{ii}| \hspace{-.5mm}<\hspace{-.5mm} w_{ii}^{k \hspace{-.2mm}-\hspace{-.2mm} 2 \hspace{-.2mm}-\hspace{-.2mm} t}(k \hspace{-.5mm}-\hspace{-.5mm} t)
\end{split}
\end{equation}
\vspace{-6mm}
\\
{for all $t+1\leq k\leq K$}.
\\
\\
Combining (\ref{eq:39}), (\ref{eq:40}), and (\ref{eq:ckt}), we have
\vspace{-2mm}
\begin{equation} \label{eq:sesvar}
    \begin{split}
    \sum_{k=0}^K\frac{\|\Delta s_{k}\|_1}{\beta_k}  
    &\leq
        2\sqrt{r}C\sum_{k=1}^{K}\sum_{t=0}^{k-1}\frac{w_{ii}^{k-1-t}}{\beta_k}\gamma_t,\\
    \sum_{k=0}^K\frac{\|\Delta x_{k}\|_1}{\beta_k} 
    & < 
    2\alpha\sqrt{r}C\sum_{k=1}^{K}\sum_{t=0}^{k-1}\frac{w_{ii}^{k-2-t}(k-t)}{\beta_k}\gamma_t.
    \end{split}
\end{equation}
\\
When $K\rightarrow \infty$ and $p\geq0$, we have
\vspace{-2mm}
\begin{equation*}
    \begin{split}
        & \sum_{k=1}^{\infty}\sum_{t=0}^{k-1}\frac{w_{ii}^{k-1-t}}{\beta_k}\gamma_t\\
    & = \gamma \sum_{k=1}^{\infty}\frac{(m+k)^q}{(m+k)^p}\sum_{t=0}^{k-1}{w_{ii}^{k-1-t}}\left(1+ \frac{k-t}{m+t}\right)^p\\
    & < \gamma \sum_{k=1}^{\infty}{(m+k)^{q-p}}\sum_{t=0}^{k-1}{w_{ii}^{k-1-t}}\left(1+ \frac{k-t}{m}\right)^p\\
    &= \gamma \sum_{k=1}^{\infty}{(m+k)^{q-p}}\sum_{t=0}^{k-1}{w_{ii}^t}\left(1+ \frac{t+1}{m}\right)^p\\
    & = \frac{\gamma}{m^pw_{ii}^{m+1}}\hspace{-1mm}\sum_{k=m+1}^{\infty}\hspace{-2mm}{k^{q-p}}\sum_{t=0}^{k-1}{w_{ii}^{m+t+1}}(m+ t+1)^p
    \end{split}
\end{equation*}
\begin{equation} \label{eq:sensbd}
\vspace{-1mm}
\begin{split}
    & < \frac{\gamma}{m^p w_{ii}^{m+1}}\hspace{-1mm}\sum_{k=m+1}^{\infty}\hspace{-2mm}{k^{q-p}}\sum_{t=0}^{k-1}{w_{ii}^{m+t+1}}(m+ t+1)^{{\lceil p \rceil}}\\
    & < \frac{\gamma}{m^p w_{ii}^{m+1}}\hspace{-1mm}\sum_{k=m+1}^{\infty}\hspace{-2mm}{k^{q-p}}\sum_{t=1}^\infty {t^{{\lceil p \rceil}}}{w_{ii}^t},
\end{split}
\end{equation}
where ${\lceil p \rceil}$ is the least integer that is no less than $p$. Note that, the term $\sum_{t=1}^\infty{t^{{\lceil p \rceil}}}{w_{ii}^t}$ is called staircase series in \cite{edgar2018staircase} and always satisfies $\sum_{t=1}^\infty{t^{{\lceil p \rceil}}}{w_{ii}^t} = \frac{w_{ii}{\rm P}_{{\lceil p \rceil}}(w_{ii})}{(1-w_{ii})^{{\lceil p \rceil}+1}}$ for $0<w_{ii}<1$, with ${\rm P}_{{\lceil p \rceil}}(w_{ii})$ the Eulerian polynomial given in (\ref{eq:pn}) \cite{edgar2018staircase}. Hence,
(\ref{eq:sensbd}) implies
\vspace{-2mm}
\begin{equation} \label{eq:infbd1}
\sum_{k=1}^{\infty}\sum_{t=0}^{k-1}\frac{w_{ii}^{k-1-t}}{\beta_k}\gamma_t
     < \frac{\gamma{\rm P}_{{\lceil p \rceil}}(w_{ii})}{m^p w_{ii}^{m}(1-w_{ii})^{{\lceil p \rceil}+1}}\sum_{k=m+1}^{\infty}\frac{1}{k^{p-q}},
\end{equation}
the right-hand side of which is always bounded under fixed $p$ and $q$ satisfying $p-q>1$ (note that the Eulerian polynomial ${\rm P}_{{\lceil p \rceil}}(w_{ii})$ is always bounded under a fixed $p$).
\\
\\
Similarly, we have 
\vspace{-2mm}
\begin{equation} \label{eq:infbd2}
    \begin{split}
        &\sum_{k=1}^{\infty}\hspace{-.5mm}\sum_{t=0}^{k-1}\hspace{-.8mm}\frac{w_{ii}^{k-2-t}(k \hspace{-.5mm}-\hspace{-.5mm} t)}{\beta_k}\gamma_t \hspace{-.5mm}<\hspace{-1mm} 
        \sum_{k=1}^{\infty}\hspace{-.5mm}\sum_{t=0}^{k-1}\hspace{-.8mm}\frac{w_{ii}^{k-2-t}(m \hspace{-.5mm}+\hspace{-.5mm} k)}{\beta_k}\gamma_t\\
        & = \gamma \sum_{k=1}^{\infty}\frac{(m+k)^{q+1}}{(m+k)^p}\sum_{t=0}^{k-1}{w_{ii}^{k-2-t}}\left(\frac{m+k}{m+t}\right)^p
        \\
        & < \gamma \frac{{\rm P}_{{\lceil p \rceil}}(w_{ii})}{m^p w_{ii}^{m+1}(1-w_{ii})^{{\lceil p \rceil}+1}}\sum_{k=m+1}^{\infty}\frac{1}{k^{p-q-1}}.
    \end{split}
\end{equation}
When $p-q>2$, $\sum_{k=m+1}^{\infty}\frac{1}{k^{p-q-1}}$ is bounded, and hence the right-hand side of (\ref{eq:infbd2}) is always bounded under the corollary conditions. 
\\
\\
Therefore, when the number of iterations $K\rightarrow\infty$, under
$p\geq0$ and $q<p-2$, condition (\ref{eq:thmDPinf}) of $b_\eta$ and $b_\xi$ ensures Proposition \ref{prop:DP} and further $\epsilon$-differential privacy of Algorithm \ref{alg:DP}.
\end{proof}

%% file: 4_opt.tex
\section{Influence of coupling weights on optimization accuracy} \label{sec:weight}
\vspace{-1mm}
{In the case of a constant stepsize, the privacy budget tends to infinity when the number of iterations tends to infinity. However, since the algorithm has a linear (exponential) convergence rate (as proven later in Theorem~\ref{thmConv0p}), we can first determine the number of iterations based on a given specification of the optimization error, and then calculate the privacy budget in the $K$ iterations. }

\subsection{Convergence analysis for constant stepsize and noise}
\begin{theorem} \label{thmConv0p} 
{Under Assumptions 1, 3, and 4,
Algorithm 1 with $\gamma_k = 1$ and $\beta_k = 1$ converges at a linear rate $\mathcal{O}(\rho({A})^k)$ if the stepsize $\alpha$ satisfies 
\vspace{-3mm}
\begin{equation*}
    \alpha<\min\left\{\frac{1}{\mu+L},\frac{1-\rho_w^2}{4\sqrt{2}d_IL},\sqrt{\frac{2d_3}{d_2+\sqrt{d_2^2+4d_1d_3}}}\right\},
\end{equation*}
where $d_1 = \frac{128 L^6d_I^2}{\mu(1-\rho_w^2)^2}$, $ d_2 = {8\mu L^2d_I^2}+ \frac{128\mu d_I^2L^2}{(1-\rho_w^2)^2}$, and $d_3 = \frac{\mu(1-\rho_w^2)^2}{4}$.}
{Moreover, we have
\vspace{-3mm}
\begin{equation*}
\begin{split}
    &\lim_{k\rightarrow\infty}\hspace{-1mm}\left\{\sup_{l\geq k}\mathbb{E}[\|\bar{x}_l\hspace{-.5mm}-\hspace{-.5mm}x^\ast\|_2^2]\hspace{-.5mm}\right\}\hspace{-.5mm}=\hspace{-1mm}\limsup_{k\rightarrow \infty} \mathbb{E}[\|{\bar{x}_k \hspace{-.5mm}-\hspace{-.5mm} {x^\ast}}\|_2^2]\hspace{-.5mm}\leq\hspace{-.5mm}\theta_1,\\
    &\lim_{k\rightarrow\infty}\hspace{-1mm}\left\{\sup_{l\geq k}\mathbb{E}[\|\boldsymbol{x}_l\hspace{-.5mm}-\hspace{-.5mm}\boldsymbol{1}\bar{x}_l\|^2_2]\hspace{-.5mm}\right\}\hspace{-.5mm}=\hspace{-1mm}\limsup_{k\rightarrow \infty} \mathbb{E}[\|{\boldsymbol{x_k} \hspace{-.5mm}-\hspace{-.5mm} \boldsymbol{1}{\bar{x}_k}}\|_2^2]\hspace{-.5mm}\leq\hspace{-.5mm}\theta_2,
\end{split}
\end{equation*}
where $\theta_1$ and $\theta_2$ are the first and second elements of the vector $(I-A)^{-1}B$, respectively, with $A$ and $B$ given in (\ref{eq:AB}) and $\theta_1$ and $\theta_2$ given in (\ref{eq:theta1}).}
\end{theorem}
\vspace{-2mm}
\begin{proof}
Under Assumption \ref{ass:network}, the following relationships hold:
\vspace{-2mm}
\begin{equation*} 
    \|W_o\|_2^2 = \sum_i |\lambda_i^o|^2\leq n\rho(W_o)^2,\quad \|\boldsymbol{1}^TW_o\|_2^2\leq n^2\rho(W_o)^2.
    \vspace{-3mm}
\end{equation*}
{Following Lemma~\ref{lemmaIneq}, when $\gamma_k =\beta_k = 1$, and $\alpha<\frac{1}{L+\mu}$, we have $\tilde{\alpha}_k = \alpha$, $\alpha\mu<1$, $\alpha L <1$, and further
\vspace{-2mm}
    \begin{equation} \label{eq:kstepbd}
        \left[  \begin{matrix} \vspace{-1mm}
        \mathbb{E}[\|\bar{x}_{k}\hspace{-0.5mm}-x^\ast\|_2^2]\\
        \mathbb{E}[\|\boldsymbol{x}_{k}\hspace{-0.5mm}-\hspace{-0.5mm}\boldsymbol{1}\bar{x}_{k}\|_2^2]\\
        \mathbb{E}[\|\boldsymbol{y}_{k}\hspace{-0.5mm}-\hspace{-0.5mm}\boldsymbol{1}\bar{y}_{k}\|_2^2]
        \end{matrix}\right]\hspace{-1mm}\leq
        \hspace{-0.5mm}A^k\hspace{-1.5mm}\left[  \begin{matrix}
        \mathbb{E}[\|\bar{x}_{0}\hspace{-0.5mm}-\hspace{-0.5mm}x^\ast\|_2^2]\\
        \mathbb{E}[\|\boldsymbol{x}_{0}\hspace{-0.5mm}-\hspace{-0.5mm}\boldsymbol{1}\bar{x}_{0}\|_2^2]\\
        \mathbb{E}[\|\boldsymbol{y}_{0}\hspace{-0.5mm}-\hspace{-0.5mm}\boldsymbol{1}\bar{y}_{0}\|_2^2]
        \end{matrix}\right]\hspace{-1.5mm} +\hspace{-1.5mm} \sum_{l=0}^{k-1} A^lB,
        \vspace{-2mm}
    \end{equation}
    with
    \vspace{-4mm}
    \begin{equation} \label{eq:AB}
    \begin{split}
        A & =
        \left[ \begin{matrix}
        1-{\alpha}\mu & \frac{2\alpha L^2}{\mu n} & 0\\
        0 & \frac{1+\rho_w^2}{2} & \frac{2\alpha^2}{1-\rho_w^2}\\
         \frac{32n{\alpha}^2 L^4d_I^2}{1-\rho_w^2} & \frac{64d_I^2L^2}{1-\rho_w^2} & \frac{1+\rho_w^2}{2} +  \frac{16{\alpha}^2L^2d_I^2}{1-\rho_w^2}
        \end{matrix}\right] 
        ,\\
        B & = \left[ \begin{matrix}
        \alpha^2\sigma_\eta^2 + \sigma_\xi^2 \\
        nd_I^2\sigma_\xi^2 \\
        \frac{24nd_I^2}{1-\rho_w^2}\sigma_\eta^2 + \frac{4nL^2d_I^2}{1-\rho_w^2}\sigma_\xi^2
        \end{matrix}\right]\rho(W_o)^2.
    \end{split}
    \vspace{-6mm}
    \end{equation}}
\\{Therefore, $\|\bar{x}_k - x^\ast\|_2^2$ and $\|\boldsymbol{x}_k-\boldsymbol{1}\bar{x}_k\|_2^2$ converge linearly to a neighborhood of zero bounded by $\theta_1$ and $\theta_2$ with the rate $\mathcal{O}(\rho(A)^k)$ if $\rho(A)<1$. When $\alpha<\frac{1}{L+\mu}<\frac{1}{\mu}$, A is nonnegative and irreducible}.   {Denote $A_{ij}$ as the element of $A$ on the $i^{th}$ row and $j^{th}$ column.} According to Lemma 3 in \cite{pu2020robust}, A necessary and sufficient condition for $\rho(A)<1$ is the elements of $A$ satisfying
$A_{11}<1$, $A_{22}<1$, $A_{33}<1$, and $\det(I-A)>0$. The inequalities
{ 
\vspace{-3mm}
\begin{equation} \label{eq:alphaBD1}
    \alpha< \frac{1}{\mu+L}, \ \ \alpha<\frac{1-\rho_w^2}{4\sqrt{2}Ld_I}
    \vspace{-2mm}
\end{equation} guarantee $A_{11}<1$ and $A_{33}<1$ respectively. In addition, $\rho_w<1$ directly inferred from Assumption~\ref{ass:network} ensures $A_{22}<1$. Also $\det(I-A)=-\alpha(d_1\alpha^4+d_2\alpha^2-d_3)>0$}
{
is guaranteed by
\vspace{-4mm}
\begin{equation} \label{eq:alphaBD2}
    \alpha^2\leq\frac{2d_3}{d_2+\sqrt{d_2^2+4d_1d_3}},
    \vspace{-2mm}
\end{equation}
where $d_1 = \frac{128 L^6d_I^2}{\mu(1-\rho_w^2)^2}$, $d_2 = {8\mu L^2d_I^2}+ \frac{128\mu d_I^2L^2}{(1-\rho_w^2)^2}$, and \\
$d_3 = \frac{\mu(1-\rho_w^2)^2}{4}$.
(\ref{eq:alphaBD1}) and (\ref{eq:alphaBD2}) give the bound on stepsize to ensure $\rho(A)<1$ and further linear convergence}.

By direct expansion, we have
\vspace{-2mm}
\begin{equation} \label{eq:theta1}
\vspace{-3mm}
\begin{split}
    [(I\hspace{-.5mm}-\hspace{-.5mm}A)^{-1}\hspace{-.5mm}B]_1 
    \hspace{-.8mm}&=\hspace{-.8mm} \frac{\alpha_1 \hspace{-.3mm}T_w^4 \hspace{-.7mm}+\hspace{-.5mm} \alpha_2 T_w^3 \hspace{-.7mm}+\hspace{-.5mm} \alpha_3 T_w^2 \hspace{-.7mm}+\hspace{-.5mm} \alpha_4 T_w \hspace{-.7mm}+\hspace{-.5mm} \alpha_5}{c_1 T_w^4 \hspace{-.5mm}+\hspace{-.5mm} c_2 T_w^2 \hspace{-.5mm}+\hspace{-.5mm} c_3}\rho({W_o})^2,\\
    [(I\hspace{-.5mm}-\hspace{-.5mm}A)^{-1}\hspace{-.5mm}B]_2
    \hspace{-.5mm}&=\hspace{-.5mm} \frac{\beta_1 T_w^3 \hspace{-.5mm}+\hspace{-.5mm} \beta_2 T_w^2 \hspace{-.5mm}+\hspace{-.5mm} \beta_3 T_w \hspace{-.5mm}+\hspace{-.5mm} \beta_4}{c_1 T_w^4 \hspace{-.5mm}+\hspace{-.5mm} c_2 T_w^2 \hspace{-.5mm}+\hspace{-.5mm} c_3}\rho({W_o})^2,
\end{split}
\vspace{-6mm}
\end{equation}
where $T_w \hspace{-.5mm}\triangleq\hspace{-.5mm} 1 - \rho_w^2$,
$\alpha_1 \hspace{-.5mm}=\hspace{-.5mm}  \frac{\alpha^2}{4}\sigma_\eta^2 + \frac{1}{4}\sigma_\xi^2$, 
$\alpha_2 \hspace{-.5mm}=\hspace{-.5mm} \frac{\alpha{d_I^2}L^2}{\mu}\sigma_\xi^2$, \\
$\alpha_3 \hspace{-.5mm}=\hspace{-.5mm} -8{d_I^2}\alpha^4L^2\sigma_\eta^2 \hspace{-.5mm}-\hspace{-.5mm}8{d_I^2}\alpha^2L^2\sigma_\xi^2$, 
$\alpha_4 \hspace{-.5mm}=\hspace{-.5mm}  -\frac{32\alpha^3 {d_I^4}L^4}{\mu}\sigma_\xi^2$, 
$\alpha_5 \hspace{-.5mm}=\hspace{-1mm} \left(\hspace{-.5mm}\frac{96\alpha^3\hspace{-.5mm}{d_I^2}L^2}{\mu} \hspace{-.5mm}-\hspace{-.5mm}  128 \alpha^4 {d_I^2}L^2\hspace{-.5mm}\right)\hspace{-.5mm}\sigma_\eta^2 \hspace{-.5mm}+\hspace{-.5mm} \left(\hspace{-.5mm}\frac{16\alpha^3\hspace{-.5mm}{d_I^2}L^4}{\mu}\hspace{-.5mm}-\hspace{-.5mm} 128 \alpha^2 {d_I^2}L^2\hspace{-.5mm}\right)$
$\times \sigma_\xi^2$,
$\beta_1 \hspace{-.5mm}=\hspace{-.5mm} \frac{\alpha\mu}{2}n{d_I^2}\sigma_\xi^2$, 
$\beta_2 \hspace{-.5mm}=\hspace{-.5mm} 0$, 
$\beta_3 \hspace{-.5mm}= \hspace{-.5mm}-16n\alpha^3\mu{d_I^4}L^2\sigma_\xi^2$, \\
$\beta_4 \hspace{-1mm}=\hspace{-1.2mm} \left(\hspace{-.5mm}64n\alpha^6\hspace{-.5mm}{d_I^2}\hspace{-.3mm}L^4 \hspace{-1mm}+\hspace{-.5mm} 48\mu n \alpha^3\hspace{-.5mm}{d_I^2}\hspace{-.5mm}\right)\hspace{-.8mm}\sigma_\eta^2 \hspace{-.5mm}+\hspace{-.5mm} \left(\hspace{-.5mm}{64n\alpha^4\hspace{-.5mm}{d_I^2}\hspace{-.3mm}L^4} \hspace{-1mm}+\hspace{-.5mm} 8\mu n\alpha^3\hspace{-.5mm}{d_I^2}\hspace{-.3mm}L^2\hspace{-.5mm}\right)$ $\times \sigma_\xi^2$,
$c_1 \hspace{-.5mm}=\hspace{-.5mm} \frac{\alpha\mu}{4}$, 
$c_2 \hspace{-.5mm}=\hspace{-.5mm} -8\alpha^3\hspace{-.5mm}\mu{d_I^2}L^2$, and
$c_3 \hspace{-.5mm}=\hspace{-.5mm} - \frac{128\alpha^5\hspace{-.5mm}{d_I^2}L^6}{\mu} \hspace{-.5mm}$ $- 128\alpha^3\hspace{-.5mm} \mu {d_I^2}L^2$.
\end{proof}
\vspace{-5mm}
From Theorem \ref{thmConv0p}, it can be seen that under the differential-privacy design, the optimization error satisfies 
\vspace{-2mm}
\begin{equation} \label{eq:bound}
\begin{split}
    &\limsup_{k\rightarrow \infty}
      \mathbb{E}[\|{\boldsymbol{x_k}-\boldsymbol{1}{x^\ast}}\|_2^2]\\
    & \leq \limsup_{k\rightarrow \infty} \mathbb{E}[2\|{\boldsymbol{1}\bar{x}_k-\boldsymbol{1}{x^\ast}}\|_2^2 + 2\|{\boldsymbol{x_k}-\boldsymbol{1}{\bar{x}_k}}\|_2^2]\\
    & \leq 2n\theta_1+ 2\theta_2.
\end{split}
\vspace{-5mm}
\end{equation}
(\ref{eq:bound}) indicates that the final error is affected by parameters in A and B, which in turn are determined by the objective functions, inter-agent coupling weights, and stepsize $\alpha$. Inspired by this observation, we apply theoretical analysis on how the coupling weight affects the optimization error induced by the differential-privacy noises.  
\vspace{-1mm}
\subsection{Results for a general coupling weight matrix $W$}
\vspace{-1mm}
We denote the optimization error as $\theta \triangleq 2n\theta_1+2\theta_2$, where $\theta_1$ and $\theta_2$ are defined in Theorem \ref{thmConv0p}. Next, we theoretically characterize the relationship between the inter-agent coupling weight $W$ and the optimization error $\theta$:
\vspace{2mm}
\begin{theorem} \label{thmrho} 
{Suppose Assumptions 1, 3, and 4 hold and the stepsize $\alpha$ satisfies 
\vspace{-3mm}
\begin{equation*}
    \alpha<\min\left\{\frac{1}{\mu+L},\frac{1-\rho_w^2}{8d_IL},\sqrt{\frac{2d_3}{d_2+\sqrt{d_2^2+4d_1d_3}}}\right\},
    \vspace{-2mm}
\end{equation*}}\\
the differential-privacy induced optimization error $\theta$ of Algorithm 1 (under $\gamma_k=\beta_k=1$) increases monotonically with {increases} in the spectral radius $\rho_w \triangleq  \rho(W-\frac{\boldsymbol{1}\boldsymbol{1}^T}{n})$ and the spectral radius of ${W_o}$, i.e.,
\vspace{-3mm}
\begin{equation}
    \frac{\partial\theta}{\partial \rho_w}>0, \ \ \frac{\partial{\theta}}{\partial{\rho\left({W_o}\right)}}\geq0.
    \vspace{-2mm}
\end{equation} 
\end{theorem}
\vspace{-5mm}
\begin{proof}
    See Appendix C.
\end{proof}

%% file: 5_simulation.tex
\section{Numerical simulations}\label{sec:simu}

\subsection{{Comparisons with existing DP approaches for distributed optimization}\label{sec:sim3}}

In this section, we compare the proposed algorithm with existing differentially private solutions (including \cite{huang2015differentially} for the static-consensus based distributed optimization and \cite{ding2021differentially}  for the gradient-tracking based distributed optimization). 

We use the optimal rendezvous problem (example 1 in \cite{huang2015differentially}), where the local objective function is given as $f_i(x_i)\triangleq \|x_i - a\|^2$, with $a\in\mathbb{R}^2$ denoting the position of an assembly point.

\begin{figure}
\begin{center}
\includegraphics[width=7.2cm]{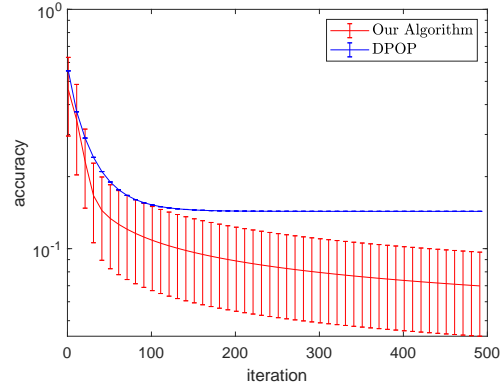}    
\caption{Comparison of our Algorithm~\ref{alg:DP} with DPOP in \cite{huang2015differentially} under the same privacy budget.}  
\label{fig:DP}                                 
\end{center}                                 
\end{figure}

{To achieve rigorous $\epsilon$-differential privacy, the approach in \cite{huang2015differentially} (hereafter referred to as DPOP) requires the stepsize to decay geometrically. In the comparison, we fix the number of iterations to $K=500$ and use the same coupling weights and privacy budget for both our algorithm and PDOP. We set the stepsize of our approach as $\alpha=0.06$, $\gamma_k = \frac{2}{(1+k)^{1.1}}$, and $\beta_k = \frac{1}{(1+k)^{0.05}}$}, and select noise variances in such a way that both algorithms have the same privacy budget. The comparison results of $500$ runs in Fig.~\ref{fig:DP} show that our algorithm yields a much better optimization accuracy under the same privacy budget.

\begin{figure}
\begin{center}
\includegraphics[width=7.2cm]{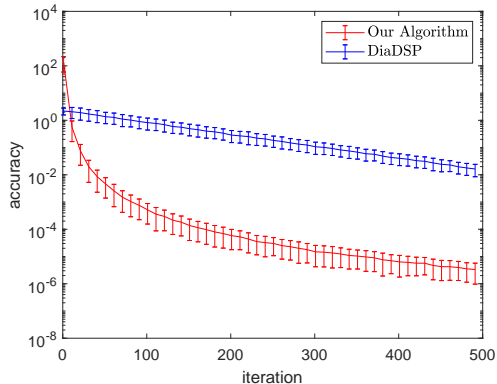}    
\caption{Comparison of our Algorithm~\ref{alg:DP} with DiaDSP in \cite{ding2021differentially} under the same privacy budget.}  
\label{fig:DP2}                                 
\end{center}                                 
\end{figure}

{
Under the same privacy budget, we also compare our algorithm with the approach in \cite{ding2021differentially} (hereafter referred to as DiaDSP). Note that \cite{ding2021differentially} restricts adjacent objective functions to have identical Lipschitz constants and convexity parameters to allow the stepsize to be constant. For a fair comparison, we use the constant stepsize $0.1$ for both algorithms. The noise decaying factor in our algorithm is set as $\beta_k=\frac{1}{(1+k)^{0.8}}$.  We fix the number of iterations to $K=500$ and select noise variances in such a way that both algorithms have the same privacy budget. The simulation results of {$500$} runs are given in Fig.~\ref{fig:DP2}. It is clear that the proposed algorithm has a much better optimization accuracy under the same privacy budget.  }

\subsection{Influence of $\rho_w$ and $\rho(W_o)$ on accuracy}\label{sec:sim1}

We use the ridge regression problem to evaluate our theoretical results for the influence of $\rho_w$ and $\rho(W_o)$ on optimization accuracy in Sec.~\ref{sec:weight}. In the distributed ridge regression problem, the objective function of the $i^{th}$ agent is given by 
$f_i\left(x\right) = \left(u_i^Tx - v_i\right)^2 + \rho\|x\|^2$ \cite{pu2020robust},
where $u_i\in\mathbb{R}^r,v_i\in\mathbb{R}$ are agent $i$'s features and observed output, respectively. Note that here with a slight abuse of notation, we use $\rho>0$ to represent a penalty parameter. For each agent, $u_i\in[-1,1]^r$ and $\Tilde{x}_i\in[0,10]^r$ are randomly generated from the uniform distribution. And in the simulation, we set $v_i=u_i^T\Tilde{x}_i+\zeta_i$, where $\zeta$ follows Gaussian distribution with mean $0$ and standard derivation $25$. This problem has a unique optimal solution
$x^\ast=\left(\sum^n_{i=1}[u_iu_i^T]+n\rho I\right)^{-1}\sum_{i=1}^n[u_iu_i^T]\Tilde{x}_i.
$

According to Theorem \ref{thmrho}, the optimization error increases with an increase in the spectral radii of ${W_o}$ and $W-\frac{\boldsymbol{1}\boldsymbol{1}^T}{n}$, i.e., $\rho(W_o)$ and $\rho_w$. In the numerical evaluation, to ensure that the $\rho(W_o)$ and $\rho_w$ can be systematically adjusted to confirm the theoretical predictions, we consider a network of four agents connected in a ring topology: 
\vspace{-3mm}
\begin{equation} \label{eq:ring}
    {W} = { \left[\begin{array}{cccc}
         \vspace{-2mm} 1-r & rd & 0 & r(1-d)  \\
         \vspace{-2mm} rd & 1-r & r(1-d) & 0  \\
         \vspace{-2mm} 0 & r(1-d) & 1-r & rd  \\
         r(1-d) & 0 & rd & 1-r
    \end{array}\right]}
    \vspace{-2mm}
\end{equation}
where $0<d<1$. This topology allows us to tune $\rho(W_o)$ and $\rho_w$ easily by changing $d$. More specifically, under the considered constraint on $W$ {and $r\leq0.5$, we have $\rho(W_o)=r$ and $\rho_w=1-r(1-d)$}. The DP-noise variances are set as $\sigma_\eta^2=\sigma_\xi^2=0.01$. 

\begin{figure}
\begin{center}
\includegraphics[width=7.2cm]{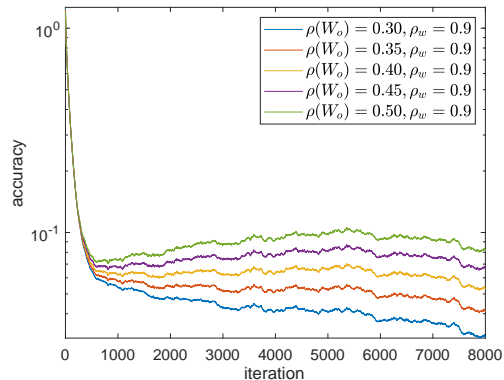}    
\caption{Influence of $\rho(W_o)$ on the optimization accuracy. $\rho_w$ is fixed to $0.9$ and $\alpha$ is fixed to $0.01$.}  
\label{fig1}                                 
\end{center}                                 
\end{figure}

\begin{figure}
\begin{center}
\includegraphics[width=7.2cm]{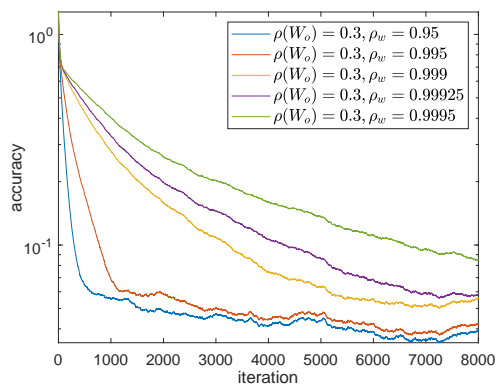}   
\caption{Influence of $\rho_w$ on optimization accuracy. $\rho(W_o)$ is fixed to $0.3$ and $\alpha$ is fixed to $0.01$.}  
\label{fig2}                                 
\end{center}                                 
\end{figure}

{The simulation results under different $\rho(W_o)$ and $\rho_w$ are shown in Fig.~\ref{fig1} and Fig.~\ref{fig2}{, respectively}. In the figures, each data point is the average of 50 runs. The results show that the optimization error indeed decreases with a decrease in $\rho(W_o)$ or $\rho_w$, corroborating the theoretical predictions in Theorem \ref{thmrho}}. 

%% file: 61_appendix1.tex
\section*{Appendix}
\section{Proof of Lemma~\ref{lemmaIneq}}
Under the assumptions $\alpha\gamma_0<\frac{2}{\mu+L}$ and $p\geq0$, we have $\Tilde{\alpha}_k<\frac{2}{\mu+L}$ with $\Tilde{\alpha}_k\triangleq\gamma_k\alpha$. Hence, for all $k\geq0$, we have
\vspace{-2mm}
\begin{equation*}
    \|\bar{x}_k-\tilde{\alpha}_k g_k-x^*\|_2\leq\left(1-\Tilde{\alpha}_k\mu\right)\|\bar{x}_k-x^*\|_2.
\end{equation*}
Under Assumption \ref{ass:network}, we have $\rho_w<1$ and $\|W\boldsymbol{w}-\boldsymbol{1}\bar{\boldsymbol{w}}\| \leq \rho_w\|\boldsymbol{w}-\boldsymbol{1}\bar{\boldsymbol{w}}\|$ for any $w\in R^{n\times r}$, where $\bar{\boldsymbol{w}}=\frac{1}{n}\boldsymbol{1}^T\boldsymbol{w}$ (See \cite{qu2017harnessing} Sec. II-B).

We denote the $\sigma$-algebra generated by $\{\boldsymbol{\eta}_0, \boldsymbol{\xi}_0,\dots,\boldsymbol{\eta}_{k-1},$ $\boldsymbol{\xi}_{k-1}\}$ as $\mathcal{F}_k$.
Similar to the proof of Lemma 2 in \cite{pu2020robust}, using the dynamics of $\bar{x}_{k+1}$ in (\ref{eq:barx}), we have
\vspace{-2mm}
\begin{equation} \label{eq:xstar}
    \begin{split}
    \vspace{-1mm}
        & \mathbb{E}[\|\bar{x}_{k+1}  -x^*\|_2^2|\mathcal{F}_k] \\
        & \hspace{-1mm}\leq (1-\Tilde{\alpha}_k\mu)\mathbb{E}[\|\bar{x}_k \hspace{-.5mm}-\hspace{-.5mm} x^*\|_2^2|\mathcal{F}_k] \hspace{-.5mm}+\hspace{-.5mm} \frac{\alpha^2\beta_k^2}{n^2}\|v\|_2^2\sigma_\eta^2 \\
        & \hspace{2mm} \hspace{-.5mm}+\hspace{-.8mm} \frac{\Tilde{\alpha}_k L^2\hspace{-.5mm}(1\hspace{-.8mm}+\hspace{-.5mm}\Tilde{\alpha}_k\mu)}{\mu n}\mathbb{E}\hspace{-.5mm}[\|\boldsymbol{x}_k \hspace{-.8mm}-\hspace{-.8mm} \boldsymbol{1}\bar{x}_k\|_2^2|\mathcal{F}_k]\hspace{-.8mm}+\hspace{-.8mm} \frac{\beta_k^2}{n^2}\hspace{-.5mm}\|v\|^2_2\sigma_\xi^2.
    \end{split}
    \vspace{-2mm}
\end{equation}
Since $\mathbb{E}[\|\left(I-\frac{\boldsymbol{1}\boldsymbol{1}^T}{n}\right){W_o}\boldsymbol{\xi}_k\|_2^2|\mathcal{F}_k]\leq d_I^2\|{W_o}\|_2^2\sigma_\xi^2$, the dynamics of $\boldsymbol{x}_{k+1}-\boldsymbol{1}\bar{x}_{k+1}$ and $\boldsymbol{y}_{k+1}-\boldsymbol{1}\bar{y}_{k+1}$ in (\ref{eq:dynx}) and (\ref{eq:dyny}) satisfy
\vspace{-2mm}
\begin{equation} \label{eq:xbar}
    \begin{split}
        & \mathbb{E}[\|\boldsymbol{x}_{k+1} 
        -\boldsymbol{1}\bar{x}_{k+1}\|_2^2|\mathcal{F}_k]\\
        & \leq \mathbb{E}[\frac{1+\rho_w^2}{2}\|\boldsymbol{x}_k-\boldsymbol{1}\bar{x}_k\|_2^2  \hspace{-.5mm}+\hspace{-.5mm}\frac{\alpha^2\left(1+\rho_w^2\right)}{1-\rho_w^2}\|\boldsymbol{y}_k-\boldsymbol{1}\bar{y}_k\|_2^2|\mathcal{F}_k]\\
        &\hspace{3.5mm} + \beta_k^2 d_I^2\|{W_o}\|_2^2\sigma_\xi^2,
    \end{split}
    \vspace{-8mm}
\end{equation}
and
\vspace{-4mm}
\begin{equation} \label{eq:bddyny}
    \begin{split}
    \vspace{-1mm}
        \|\boldsymbol{y}_{k+1} 
        \hspace{-.5mm}-\hspace{-.5mm} \boldsymbol{1}\bar{y}_{k+1}\|_2^2 & \leq \frac{1+\rho_w^2}{2}\|\boldsymbol{y}_k-\boldsymbol{1}\bar{y}_k\|_2^2  \\
        & \hspace{3.5mm}+\hspace{-.5mm} \frac{1\hspace{-.5mm}+\hspace{-.5mm}\rho_w^2}{1\hspace{-.5mm}-\hspace{-.5mm}\rho_w^2}d_I^2\|\Tilde{\nabla}F(\boldsymbol{x}_{k+1}) \hspace{-.5mm}-\hspace{-.5mm} \Tilde{\nabla}F(\boldsymbol{x}_k)\|_2^2. 
        \vspace{-4mm}
    \end{split}
\end{equation}
Denoting $\Tilde{\nabla}_k \triangleq \Tilde{\nabla} F(\boldsymbol{x}_k)$ and $\nabla_k \triangleq \nabla F(\boldsymbol{x}_k)$, we have
\vspace{-2mm}
\begin{equation} \label{eq:graddiff}
\begin{split}
    &\mathbb{E}[\|\Tilde{\nabla}_{k+1} \hspace{-.5mm}-\hspace{-.5mm} \Tilde{\nabla}_k\|_2^2|\mathcal{F}_k]
    \\
    & = \mathbb{E}[\|\gamma_{k+1}\nabla_{k+1} \hspace{-.5mm}-\hspace{-.5mm}\gamma_k\nabla_k\|_2^2|\mathcal{F}_k] \hspace{-.5mm}+\hspace{-.5mm} \|W_0\|^2_2(\beta^2_{k+1}\hspace{-.5mm}+\hspace{-.5mm}\beta_k^2)\sigma_\eta^2
    \\ 
    & \hspace{-.5mm}+\hspace{-.5mm} 2 \mathbb{E}[\langle\gamma_{k+1}\nabla_{k+1}\hspace{-.5mm}-\hspace{-.5mm}\gamma_k\nabla_k,W_0\beta_{k+1}\eta_{k+1}\hspace{-.5mm}-\hspace{-.5mm}W_0\beta_k\eta_{k}\rangle|\mathcal{F}_k]\\
    & \hspace{-.5mm}\leq \hspace{-.5mm} 2\gamma_{k+1}^2
        \mathbb{E}[\|\nabla_{k+1} \hspace{-.5mm}-\hspace{-1mm} \nabla_k\|_2^2|\mathcal{F}_k] 
        \hspace{-.5mm}+\hspace{-.5mm} 2(\gamma_k\hspace{-.5mm}-\hspace{-.5mm}\gamma_{k+1})^2 \mathbb{E}[\|\nabla_k\|^2_2|\mathcal{F}_k]\\
    & + 2\gamma_{k+1}\mathbb{E}[\langle\nabla_{k+1}, - W_0\beta_k\eta_{k}\rangle|\mathcal{F}_k] +2\|{W_o}\|_2^2\beta_k^2\sigma_\eta^2.
\end{split}
\vspace{-2mm}
\end{equation}
Using Assumption~\ref{ass:muL} and Assumption~\ref{ass:noise}, 
the first term on the right-hand side of (\ref{eq:graddiff}) satisfies
\vspace{-2mm}
\begin{equation} \label{eq:graddiff1}
    \begin{split}
        &\mathbb{E}[\|\nabla_{k+1}  -\nabla_k\|_2^2|\mathcal{F}_k] \\
        & \leq L^2\mathbb{E}[\|\boldsymbol{x}_{k+1}-
        \boldsymbol{x}_k\|_2^2|\mathcal{F}_k]\\
        & \leq 4L^2(1+\rho_w^2) \mathbb{E}[\|\boldsymbol{x}_k -\boldsymbol{1}\bar{x}_k\|_2^2|\mathcal{F}_k]+ L^2\|{W_o}\|^2_2\beta_k^2\sigma_\xi^2\\
        & \hspace{3.5mm} +
        4\alpha^2L^2\mathbb{E}[\|\boldsymbol{y}_k -\boldsymbol{1}\bar{y}_k\|_2^2|\mathcal{F}_k]
        + 4\alpha^2L^2n \mathbb{E}[ \|\bar{y}_k\|_2^2|\mathcal{F}_k].
    \end{split}
    \vspace{-2mm}
\end{equation}
{Denote $\nabla^\star \triangleq  [\nabla f_1(x^\star), \dots, \nabla f_n(x^\star)]^T$. Based on Assumption \ref{ass:muL}, the second term on the right-hand side of (\ref{eq:graddiff}) satisfies
\vspace{-2mm}
\begin{equation} \label{eq:graddiff2}
\begin{split}
    \mathbb{E}[\|\nabla_k\|_2^2|\mathcal{F}_k] & = \mathbb{E}[\|\nabla_k-\nabla^\star + \nabla^\star\|_2^2|\mathcal{F}_k]\\
    & \leq
    2L^2\mathbb{E}[\|\boldsymbol{x}_{k}- 
    \boldsymbol{1}[x^\star]^T\|_2^2|\mathcal{F}_k] + 2c^\star\\
    & \leq 
    4L^2\mathbb{E}[\|\boldsymbol{x}_{k}- \boldsymbol{1}\bar{x}_k\|^2_2|\mathcal{F}_k] + 2c^\star\\
    & \hspace{3.5mm} + 4nL^2\mathbb{E}[\|\bar{x}_k - 
    {x}^\star\|_2^2|\mathcal{F}_k] ,
    \vspace{-2mm}
\end{split}
\end{equation}
where $c^\star \triangleq \sum_{i=1}^n\|f_i(x^\star)\|_2^2$, which is a constant for any fixed Problem \ref{eq:DistOpt1} under Assumption \ref{ass:muL}.}

Following the derivations in Lemma 8 of \cite{pu2021distributed}, 
we can bound the third term on the right-hand side of (\ref{eq:graddiff}) as follows:
\vspace{-2mm}
\begin{equation} \label{eq:graddiff3}
\begin{split}
    & \mathbb{E}[\langle\nabla_{k+1}, - W_0\beta_k\eta_{k}\rangle|\mathcal{F}_k]\\
    & \leq \hspace{-.8mm} \alpha \hspace{-1mm}  \sum_{i=1}^n \mathbb{E}\hspace{-1.5mm}\left.\left[ \hspace{-.5mm}L\left\|\sum_{j=1,j\neq i}^p \hspace{-2mm}w_{ij}\beta_k\eta_{j,k}\right\|^2\right|\mathcal{F}_k\right] \hspace{-1mm}\leq\hspace{-.5mm} \alpha L \beta_k^2 \|W_0\|^2_2\sigma_\eta^2.
\end{split}
\vspace{-2mm}
\end{equation}
Combining (\ref{eq:graddiff}), (\ref{eq:graddiff1}), (\ref{eq:graddiff2}), and (\ref{eq:graddiff3}) with the following inequality:
\vspace{-2mm}
\begin{equation*}
    \begin{split}
        &\mathbb{E}[\|\bar{y}_k\|_2^2|\mathcal{F}_k]
         \leq \frac{\|v\|_2^2}{n^2}\sigma_\eta^2 + \frac{2L^2}{n}\mathbb{E}[\|\boldsymbol{x}_k - \boldsymbol{1}\bar{x}_k\|_2^2|\mathcal{F}_k] \\
        & \hspace{2.2cm} + 2L^2\mathbb{E}[\|\bar{x}_k-x^*\|_2^2|\mathcal{F}_k],
    \end{split}
\end{equation*} we have 
\vspace{-2mm}
\begin{equation*}
    \begin{split}
        & \mathbb{E}[\|\Tilde{\nabla}_{k+1}-\Tilde{\nabla}_k\|_2^2|\mathcal{F}_k] \\
        & \leq \left(16\gamma_{k+1}^2 \Tilde{\alpha}_k^2L^4n + 
        {8nL^2(\gamma_k-\gamma_{k+1})^2}\right)
        \mathbb{E}[\|\bar{x}_k \hspace{-.5mm}-\hspace{-.5mm} x^*\|_2^2|\mathcal{F}_k]\\
        & \hspace{3.5mm} +\hspace{-.5mm} \left(2\gamma_{k+1}^2 (4L^2(1+\rho_w^2)+8\Tilde{\alpha}_k^2L^4) + {8L^2(\gamma_k-\gamma_{k+1})^2}\right)\\
        & \hspace{3.5mm}
        \times\mathbb{E}[\|\boldsymbol{x}_k -\boldsymbol{1}\bar{x}_k\|_2^2|\mathcal{F}_k]
        \hspace{-.5mm}+\hspace{-.5mm} 2\gamma_{k+1}^2\beta_k^2L^2\|{W_o}\|_2^2\sigma_\xi^2\\
        & \hspace{3.5mm}
        + \hspace{-.5mm} 8\Tilde{\alpha}_{k+1}^2L^2\mathbb{E}[\|\boldsymbol{y}_k -\boldsymbol{1}\bar{y}_k\|_2^2|\mathcal{F}_k]
        + {4(\gamma_k-\gamma_{k+1})^2c^*}\\ 
        & \hspace{3.5mm}
        + \hspace{-1mm} \left(\hspace{-.5mm}\frac{8\Tilde{\alpha}_{k+1}^2\beta_k^2L^2\|v\|_2^2}{n}
        \hspace{-.5mm}+\hspace{-.5mm}  2\Tilde{\alpha}_{k+1}\beta_k^2 L\|W_0\|^2_2 \hspace{-.5mm}+\hspace{-.5mm} 2\beta_k^2\|{W_o}\|_2^2\hspace{-.5mm}\right)\hspace{-1mm}\sigma_\eta^2.
    \end{split}
    \vspace{-2mm}
\end{equation*}
Further using the relationship $1+\rho_w^2<2$, we have the following relationship based on (\ref{eq:bddyny}):
\vspace{-2mm}
\begin{equation} \label{eq:ybar}
    \begin{split}
        & \mathbb{E}[\|\boldsymbol{y}_{k+1}  -\boldsymbol{1}\bar{y}_{k+1}\|_2^2|\mathcal{F}_k]\\ 
        & \leq \left(\frac{32n\gamma_{k+1}^2\Tilde{\alpha}_k^2L^4d_I^2}{1-\rho_w^2} + {\frac{16(\gamma_k-\gamma_{k+1})^2nL^2d_I^2}{1-\rho_w^2}}\right)\\
        & \hspace{3.5mm} \times
        \mathbb{E}[\|\bar{x}_k - x^*\|_2^2|\mathcal{F}_k] +  \left(\frac{32\gamma_{k+1}^2d_I^2(L^2+\Tilde{\alpha}_k^2L^4)}{1-\rho_w^2} \right.\\
        & \hspace{3.5mm} \left. + {\frac{16(\gamma_k-\gamma_{k+1})^2L^2d_I^2}{1-\rho_w^2}}\right)
        \mathbb{E}[\|\boldsymbol{x}_k -\boldsymbol{1}\bar{x}_k\|_2^2|\mathcal{F}_k]\\
        & \hspace{3.5mm} + \left(\frac{1+\rho_w^2}{2}+\frac{16\Tilde{\alpha}_{k+1}^2L^2d_I^2}{1-\rho_w^2}\right)\mathbb{E}[\|\boldsymbol{y}_k -\boldsymbol{1}\bar{y}_k\|_2^2|\mathcal{F}_k]\\
        & \hspace{3.5mm} + \frac{4\beta_k^2d_I^2\sigma_\eta^2}{1-\rho_w^2}\left(\frac{4\Tilde{\alpha}_{k+1}^2L^2\|v\|_2^2}{n}
        +  (1 + \Tilde{\alpha}_{k+1} L )\|{W_o}\|_2^2\right)\\
        & \hspace{3.5mm} + \frac{4\gamma_{k+1}^2\beta_k^2L^2d_I^2}{1-\rho_w^2}\|{W_o}\|_2^2\sigma_\xi^2
        + {\frac{8(\gamma_k-\gamma_{k+1})^2d_I^2c^\star}{1-\rho_w^2}}.
    \end{split}
    \vspace{-2mm}
\end{equation}

The result in (\ref{eq:linear}) and (\ref{eq:A_k}) is obtained by taking the full expectations of (\ref{eq:xstar}), (\ref{eq:xbar}), and (\ref{eq:ybar}) via the law of total expectation.

%% file: 62_appendix2.tex
{\section{Proof of Theorem~\ref{thmConv}}\label{sec:proofthm2}}

Define $U_k\triangleq\mathbb{E}[\|\bar{x}_{k}\hspace{-0.5mm}-x^\ast\|_2^2]$, $X_k\triangleq\mathbb{E}[\|\boldsymbol{x}_{k}\hspace{-0.5mm}-\hspace{-0.5mm}\boldsymbol{1}\bar{x}_{k}\|_2^2]$, and $Y_k\triangleq\mathbb{E}[\|\boldsymbol{y}_{k}\hspace{-0.5mm}-\hspace{-0.5mm}\boldsymbol{1}\bar{y}_{k}\|_2^2]$.

\textbf{Step 1:} We first prove that for $p\geq 0$, the following relations hold for all $k\geq 0$:
\vspace{-2mm}
\begin{equation} \label{eq:UXYk}
    U_{k}\leq f_{k}\hat{U}, X_{k}\leq g_{k}\hat{X}, Y_{k}\leq h_{k}\hat{Y},
    \vspace{-2mm}
\end{equation}
where $f_k$, $g_k$, and $h_k$ are some positive decreasing functions of $k$, and $\hat{U}$, $\hat{X}$, and $\hat{Y}$ are some constants independent of $k$.

{
For $k = 0$, (\ref{eq:UXYk}) holds if $\hat{U}$, $\hat{X}$, and $\hat{Y}$ satisfy
\vspace{-3mm}
\begin{equation} \label{eq:UXYinit}
    \hat{U}\geq\frac{U_0}{f_0}, \hat{X}\geq\frac{X_0}{g_0}, \hat{Y}\geq\frac{Y_0}{h_0}.
    \vspace{-3mm}
\end{equation}
Given $U_{k}\leq f_{k}\hat{U}$, $X_{k}\leq g_{k}\hat{X}$, $Y_{k}\leq h_{k}\hat{Y}$, and 
\vspace{-2mm}
\begin{equation} \label{eq:alphagamma0}
\alpha\gamma_0\leq\frac{2}{\mu+L},
\vspace{-3mm}
\end{equation}
a sufficient condition ensuring $U_{k+1}\leq f_{k+1}\hat{U}$, $X_{k+1}\leq g_{k+1}\hat{X}$, and $Y_{k+1}\leq h_{k+1}\hat{Y}$ is (see Lemma 1)
\vspace{-4mm}
\begin{equation} \label{eq:UXYineq}
    \begin{split}
        U_{k+1}&\leq (1-\tilde{\alpha}_k\mu)f_k \hat{U} + \tilde{\alpha}_ka_{12}g_k\hat{X} + \beta_k^2b_1 \leq f_{k+1}\hat{U},\\
        X_{k+1}&\leq\frac{1+\rho_w^2}{2}g_k\hat{X}+a_{23}\alpha^2h_k\hat{Y}+\beta_k^2b_2\leq g_{k+1}\hat{X},\\
        Y_{k+1}& \leq\left(\gamma^2_{k+1}\tilde{\alpha}_k^2a_{31}+ (\gamma_k-\gamma_{k+1})^2a_{34}\right)f_k\hat{U}\\
        & \hspace{3.5mm} + \left(\gamma^2_{k+1}a_{32} + (\gamma_k-\gamma_{k+1})^2a_{35}\right)g_k\hat{X}\\
        & \hspace{3.5mm}+ \left(\frac{1+\rho_w^2}{2}+\tilde{\alpha}_{k+1}^2a_{33}\right)h_k\hat{Y}\\
        & \hspace{3.5mm} + \beta_k^2b_{31}+(\gamma_k - \gamma_{k+1})^2b_{32}\leq h_{k+1}\hat{Y}.
    \end{split}
    \vspace{-12mm}
\end{equation}
(\ref{eq:UXYineq}) holds if there exist $\hat{U}$, $\hat{X}$, and $\hat{Y}$ satisfying
\vspace{-2mm}
\begin{equation} \label{eq:d}
    \begin{split}
        {d}_1\hat{U}& \geq \tilde{\alpha}_ka_{12}g_k\hat{X} + \beta_k^2b_1,\\
        {d}_2\hat{X}&\geq a_{23}\alpha^2h_k\hat{Y}+\beta_k^2b_2,\\
        {d}_3\hat{Y}& \geq \left(\gamma^2_{k+1}\tilde{\alpha}_k^2a_{31} + (\gamma_k-\gamma_{k+1})^2a_{34}\right)f_k\hat{U}
        \hspace{-.5mm}\\
        & \hspace{3.5mm}+\hspace{-.5mm}\left(\gamma^2_{k+1}a_{32} + (\gamma_k-\gamma_{k+1})^2a_{35}\right)g_k\hat{X} 
        \hspace{-.5mm}\\
        & \hspace{3.5mm} +\hspace{-.5mm} \beta_k^2b_{31}
        \hspace{-.5mm}+\hspace{-.5mm}(\gamma_k \hspace{-.5mm}-\hspace{-.5mm} \gamma_{k+1})^2b_{32},
    \end{split}
\end{equation}
with
$d_1\triangleq f_{k+1} - (1-\tilde{\alpha}_k\mu)f_k$, $d_2\triangleq g_{k+1} - \frac{1+\rho_w^2}{2}g_k$, and $d_3\triangleq h_{k+1} - (\frac{1+\rho_w^2}{2}+\tilde{\alpha}_{k+1}^2a_{33})h_k$.}

{
For $p\geq 0$, $m>0$, and $k \geq 0$, $\Tilde{\alpha}_k\triangleq\alpha\gamma_k$ is bounded by $\Tilde{\alpha}_k\leq \alpha\gamma_0$.  Denoting $c_{mp}\triangleq1-\left(\frac{m}{m+1}\right)^p$,
we also have $\gamma_k - \gamma_{k+1} ={\gamma} \frac{1-\left(\frac{m+k}{m+k+1}\right)^p}{(m+k)^{p}}< {\gamma}\frac{c_{mp}}{(m+k)^{p}}$.
Given any $\tilde{d}_1$, $\tilde{d}_2$, and $\tilde{d}_3$ satisfying $0<\tilde{d}_1\leq d_1$, $0<\tilde{d}_2\leq d_2$, and $0<\tilde{d}_3\leq d_3$, 
it suffices to show that there exist $\hat{U}$, $\hat{X}$, and $\hat{Y}$ satisfying
\vspace{-3mm}
\begin{equation} \label{eq:bdineq}
    \begin{split}
    \vspace{-1mm}
    \hat{U} & \geq \left(\tilde{\alpha}_k\tilde{a}_{12}g_k\hat{X} + \beta_k^2b_1\right)/\tilde{d}_1,\\
        \hat{X} & \geq \left(a_{23}\alpha^2h_k\hat{Y}+\beta_k^2b_2\right)/\tilde{d}_2,\\
    \vspace{-1mm}
        \hat{Y} & 
        \geq 
        \hspace{-1mm}\left( \hspace{-1.5mm}\left(\gamma^2_{k+1}\tilde{\alpha}_k^2a_{31} \hspace{-1mm}+\hspace{-1mm} \frac{\gamma^2a_{34}c_{mp}^2}{(m+k)^{2p}}\hspace{-1.5mm}\right)\hspace{-1mm} f_k\hat{U} \hspace{-.5mm}+\hspace{-.5mm} \left(\gamma^2_{k+1}\tilde{a}_{32} \right.\right.\\
        & \hspace{3.5mm} \left. \left.+ \frac{\gamma^2a_{35}c_{mp}^2}{(m+k)^{2p}}\hspace{-1.5mm} \right)\hspace{-1mm}g_k\hat{X} \hspace{-.5mm}+\hspace{-.5mm} \beta_k^2\tilde{b}_{31} \hspace{-.5mm}+\hspace{-.5mm} \frac{{\gamma^2b_{32}c_{mp}^2}}{(m \hspace{-.5mm}+\hspace{-.5mm} k)^{2p}}\hspace{-1.5mm}\right)\hspace{-1mm}/\tilde{d}_3,
    \end{split}
\end{equation}
where $\tilde{a}_{12}\triangleq\frac{L^2}{\mu n}(1+\alpha\gamma_0\mu)$, $\tilde{a}_{32} \triangleq \frac{32d_I^2(L^2+\alpha^2\gamma_0^2L^4)}{1-\rho_w^2}$, and $\tilde{b}_{31}\triangleq \frac{4d_I^2\sigma_\eta^2}{1-\rho_w^2}\left(\frac{4{\alpha}^2\gamma_0^2L^2\|v\|_2^2}{n} +  (1 + {\alpha}\gamma_0 L )\|{W_o}\|_2^2\right) + \frac{4\gamma_0^2L^2d_I^2}{1-\rho_w^2}\|{W_o}\|_2^2\sigma_\xi^2$.
The exact expressions of $\tilde{d}_1$, $\tilde{d}_2$, and $\tilde{d}_3$ and the conditions ensuring the right-hand sides of (\ref{eq:bdineq}) to be bounded will be discussed later. 

Given $e_0$, $e_1$, $e_2$, $e_3$, $e_4$, $e_5$, and $e_6$ independent of $k$ such that $e_0\geq\left(\gamma^2_{k+1}\tilde{\alpha}_k^2a_{31} {+ \frac{\gamma^2a_{34}c_{mp}^2}{(m+k)^{2p}}}\right) \frac{f_k}{\tilde{d}_3}$, 
$e_1\geq\left(\gamma^2_{k+1}\tilde{a}_{32} + \frac{\gamma^2a_{35}c_{mp}^2}{(m+k)^{2p}} \right)\frac{g_k}{\tilde{d}_3}$, 
$e_2\geq\left(\beta_k^2\tilde{b}_{31} + \frac{{\gamma^2b_{32}c_{mp}^2}}{(m+k)^{2p}}\right)/\tilde{d}_3$, $e_3\geq\frac{a_{23}\alpha^2h_k}{\tilde{d}_2}$, $e_4\geq\frac{\beta_k^2b_2}{\tilde{d}_2}$, $e_5\geq\frac{\tilde{\alpha}_k\tilde{a}_{12}g_k}{\tilde{d}_1}$, and $e_6\geq\frac{\beta_k^2b_1}{\tilde{d}_1}$,
(\ref{eq:bdineq}) is guaranteed by  
\vspace{-2mm}
\begin{equation} \label{eq:eineq}
    \begin{split}
        \hat{U}& \geq e_5\hat{X} + e_6,\\
        \hat{X}& \geq e_3\hat{Y} + e_4, \\ 
        \hat{Y} & 
        \geq 
         e_0\hat{U} + e_1\hat{X} + e_2.
    \end{split}
    \vspace{-2mm}
\end{equation}
We will give the exact expressions of $e_0$ to $e_6$ under different values of $p$ later.

Combining (\ref{eq:UXYinit}) and (\ref{eq:eineq}), we can set $\hat{U}$ as
\vspace{-2mm}
\begin{equation} \label{eq:hatU}
    \hat{U} = \max \{e_5\hat{X} + e_6, U_0/f_0 \}.
\end{equation}
By plugging the value of $\hat{U}$ in (\ref{eq:hatU}) into the inequality of $\hat{X}$ and $\hat{Y}$ in (\ref{eq:eineq}), it is sufficient to show that there exists $\hat{X}$ and $\hat{Y}$ satisfying
\vspace{-3mm}
\begin{equation} \label{eq:condhatY}
    e_0\hat{U} + e_1\hat{X} + e_2\leq \hat{Y} \leq \frac{1}{e_3}(\hat{X} - e_4).
\end{equation}
When 
\vspace{-2mm}
\begin{equation} \label{eq:condtodiscuss}
    \frac{1}{e_3}-e_0e_5 - e_1>0,
\end{equation}
one can always choose $\hat{X}$ as
\vspace{-2mm}
\begin{equation} \label{eq:hatX}
\begin{split}
    \hat{X} \hspace{-1mm}=\hspace{-.5mm} & \max \hspace{-1mm}\left\{ \hspace{-1.5mm}\left(\hspace{-.5mm}\frac{1}{e_3} \hspace{-.5mm}-\hspace{-.5mm} e_0e_5 \hspace{-.5mm}-\hspace{-.5mm} e_1\hspace{-1.5mm}\right)^{\hspace{-1.5mm}-1}\hspace{-2mm}\left(\hspace{-.5mm}\frac{e_4}{e_3} + e_0e_6 +  e_2\hspace{-.5mm}\right)\hspace{-1mm}, \frac{X_0}{g_0}, \right.\\
    & \left.\hspace{-1mm}\left(\hspace{-.5mm}\frac{1}{e_3} - e_1\hspace{-.5mm}\right)^{\hspace{-1mm}-1}\hspace{-2mm}\left(\hspace{-.5mm}\frac{e_4}{e_3} + e_0\frac{U_0}{f_0} + e_2\hspace{-.5mm}\right)\hspace{-1mm}, e_3\frac{Y_0}{h_0} + e_4\hspace{-.5mm}\right\},
\end{split}
\end{equation} to satisfy (\ref{eq:UXYinit}) and 
$\frac{1}{e_3}(\hat{X} - e_4) \geq \max\{e_0\hat{U} + e_1\hat{X} + e_2, Y_0/h_0\},$
which further ensures the existence of $\hat{Y}$ satisfying (\ref{eq:UXYinit}) and (\ref{eq:condhatY}).
Then choosing $\hat{Y}$ as 
\vspace{-2mm}
\begin{equation} \label{eq:hatY}
    \hat{Y} = \frac{1}{e_3}(\hat{X}-e_4).
\end{equation}
satisfies the relationship in (\ref{eq:UXYinit}) and (\ref{eq:condhatY}).

In summary, combining (\ref{eq:UXYinit}) and (\ref{eq:eineq}), one can always find $\hat{U}$, $\hat{X}$, and $\hat{Y}$ as (\ref{eq:hatU}), (\ref{eq:hatX}), and (\ref{eq:hatY}), respectively, all of which are independent of $k$, to satisfy (\ref{eq:bdineq}) and further (\ref{eq:UXYk}) for all $k\geq 0$.

Therefore, following (\ref{eq:UXYk}), we get the error evolution of Algorithm 1 as
\vspace{-2mm}
\begin{equation} \label{eq:errorevo}
    \begin{split}
        & \mathbb{E}[\|\bar{x}_{k+1}\hspace{-0.5mm}-x^\ast\|_2^2]
        \leq f_k{\mathcal{O}(1)},\\
        & \mathbb{E}[\|\boldsymbol{x}_{k+1} -\boldsymbol{1}\bar{x}_{k+1}\|_2^2]
        \leq g_k{\mathcal{O}(1)},\\
        & \mathbb{E}[\|\boldsymbol{y}_{k+1}-\boldsymbol{1}\bar{y}_{k+1}\|_2^2]
        \leq h_k {\mathcal{O}(1)},
    \end{split}
    \vspace{-2mm}
\end{equation}
where $\mathcal{O}(1)$ is some constant independent of $k$.

{\textbf{Step 2:} In this step, we give the concrete expressions of $f_k$, $g_k$, and $h_k$ that satisfy (\ref{eq:UXYineq}) and $\tilde{d}_1$, $\tilde{d}_2$ and $\tilde{d}_3$ that satisfy (\ref{eq:bdineq}). We also give the concrete expressions of $e_0$, $e_1$, $e_2$, $e_3$, $e_4$, $e_5$, and $e_6$ satisfying (\ref{eq:eineq}). For notational convenience, we denote $m+k$ as $\ell_k$.}

{
\textbf{Case I:} For $0\hspace{-.5mm}\leq\hspace{-.5mm} p \leq 1$, we select $f_k$, $g_k$, and $h_k$ as
\vspace{-2mm}
\begin{equation} \label{eq:fgh_leq1}
    f_k=\frac{1}{(m+k)^l}, g_k=\frac{1}{(m+k)^g}, h_k=\frac{1}{(m+k)^h},
\end{equation}
where $l$, $g$, and $h$ are some positive numbers dependent on $p$ and $q$. Plugging (\ref{eq:fgh_leq1}) into (\ref{eq:d}), we get
\vspace{-2mm}
\begin{equation} \label{eq:dp0}
\begin{split}
    d_1 & 
    = \frac{\alpha\gamma\mu+\frac{{\ell_k}^{l+p}}{(\ell_{k}+1)^l}-{\ell_k}^p}{{\ell_k}^{l+p}}, 
    d_2 = \frac{2\left(\frac{\ell_{k}}{{\ell_{k}+1}}\right) ^g - (1+\rho_w^2)}{{\ell_k}^g},\\
    d_3 & = \frac{\left(\frac{\ell_{k}}{{\ell_{k}+1}}\right)^h-(\frac{1+\rho_w^2}{2}+\frac{\alpha^2\gamma^2a_{33}}{(\ell_{k}+1)^{2p}})}{{\ell_k}^h}.
\end{split}
\vspace{-10mm}
\end{equation}
1) For $p=0$, we have
\vspace{-4mm}
\begin{equation} \label{eq:d1_p1}
    \alpha\gamma\mu+\frac{{\ell_k}^l}{(\ell_{k}+1)^l}-1 > \alpha\gamma\mu + \left(\frac{m}{m+1}\right)^l -1.
    \vspace{-2mm}
\end{equation}
2) For $0<p< l\leq 1$, when ${\ell_k}\geq\frac{p}{(l-p)}>0$, we have
\vspace{-4mm}
\begin{equation} \label{eq:partiald1}
\begin{split}
    & \frac{\partial\{\frac{{\ell_k}^{l+p}}{({\ell_k}+1)^l}-{\ell_k}^p\}}{\partial {\ell_k}} \\
    & = \hspace{-.5mm}{\ell_k}^{p-1}\hspace{-1.5mm}\left(\hspace{-.8mm}\frac{{\ell_k}}{{\ell_k}\hspace{-.8mm}+\hspace{-.8mm}1}\hspace{-.5mm}\right)^{\hspace{-.8mm}l-1}\hspace{-1.8mm}\left(\hspace{-1mm}\frac{(l\hspace{-.5mm}-\hspace{-.5mm}p){\ell_k}}{({\ell_k}\hspace{-.8mm}+\hspace{-.8mm}1)^2} \hspace{-.5mm}+\hspace{-.5mm} p\frac{{\ell_k}({\ell_k}\hspace{-.8mm}+\hspace{-.8mm}2)}{({\ell_k}\hspace{-.8mm}+\hspace{-.8mm}1)^2} \hspace{-.5mm}-\hspace{-.5mm} p\hspace{-1mm}\left(\hspace{-1mm}\frac{{\ell_k}}{{\ell_k}\hspace{-.8mm}+\hspace{-.8mm}1}\hspace{-1mm}\right)^{\hspace{-1mm}1-l} \hspace{-.3mm}\right)\\
    & \geq \hspace{-.5mm}{\ell_k}^{p-1}\hspace{-1mm}\left(\frac{{\ell_k}}{{\ell_k}\hspace{-.5mm}+\hspace{-.5mm}1}\right)^{l-1}\hspace{-1.5mm}\left(\frac{(l-p){\ell_k}}{({\ell_k}\hspace{-.5mm}+\hspace{-.5mm}1)^2} \hspace{-.5mm}+\hspace{-.5mm} p\frac{{\ell_k}({\ell_k}\hspace{-.5mm}+\hspace{-.5mm}2)}{({\ell_k}\hspace{-.5mm}+\hspace{-.5mm}1)^2} \hspace{-.5mm}-\hspace{-.5mm} p \right)\\
    & = \hspace{-.5mm} {\ell_k}^{p-1}\hspace{-1mm}\left(\frac{{\ell_k}}{{\ell_k}\hspace{-.5mm}+\hspace{-.5mm}1}\right)^{l-1}\hspace{-1.5mm}\left(\frac{(l-p){\ell_k} - p}{({\ell_k}+1)^2} \right) \geq 0.
\end{split}
\vspace{-8mm}
\end{equation}
Hence, 
when $m\geq{\frac{p}{(l-p)}}$,
we have 
\vspace{-4mm}
\begin{equation} \label{eq:d1_p2}
    \alpha\gamma\mu+\frac{{\ell_k}^{l+p}}{(\ell_{k}+1)^l}-{\ell_k}^p
    \geq \alpha\gamma\mu+\frac{m^{l+p}}{(m+1)^l}-m^p.
    \vspace{-2mm}
\end{equation}
3) For $0<p\leq 1$ and $l\geq 1$, using the Bernoulli inequality $(1+x)^l\geq1+lx$ for $x\geq -1$ and $l\geq 1$, we have $\left(1-\frac{1}{{\ell_{k}+1}}\right)^l\geq 1 - l/{(\ell_{k}+1)}$, which further yields
\vspace{-4mm}
\begin{equation} \label{eq:d1_p3}
    d_1=\frac{\alpha\gamma\mu}{{\ell_k}^{l+p}}+\frac{\left(1\hspace{-.5mm}-\hspace{-.5mm}\frac{1}{{\ell_{k}+1}}\right)^l\hspace{-.5mm}-\hspace{-.5mm}1}{{\ell_k}^l}
    \geq\frac{\alpha\gamma\mu-l}{{\ell_k}^{l+p}}.
    \vspace{-2mm}
\end{equation}
Combining (\ref{eq:d1_p1}), (\ref{eq:d1_p2}), and (\ref{eq:d1_p3}), we have 
\vspace{-2mm}
\begin{equation} \label{eq:d1}
    d_1\geq \frac{\alpha\gamma\mu - c(p,l)}{{\ell_k}^{l+p}},
    \vspace{-4mm}
\end{equation} with 
\vspace{-4mm}
\begin{equation*}
    c(p,l) \hspace{-.5mm}=\hspace{-.5mm} \left\{\begin{array}{cl}
    1 \hspace{-.5mm}-\hspace{-.5mm} \left(\frac{m}{m+1}\right)^l & \text{ for } p=0, l>0\\
    m^p \hspace{-.5mm}-\hspace{-.5mm}\frac{m^{p+l}}{(m+1)^l}  & \text{ for } 0<p<l< 1, m\geq\frac{p}{l-p}\\
    l & \text{ for } 0<p\leq1\leq l
    \end{array}\right.\hspace{-.5mm}.
    \vspace{-3mm}
\end{equation*}
Combining (\ref{eq:dp0}) and (\ref{eq:d1}), for $0\leq p\leq1$ , we select $\tilde{d}_1\leq d_1$, $\tilde{d}_2\leq d_2$, and $\tilde{d}_3 \leq d_3$ in (\ref{eq:bdineq}) as 
\vspace{-3mm}
\begin{equation*}
\begin{split}
    & \Tilde{d}_1 = \frac{\alpha\gamma\mu - c(p,l)}{{\ell_k}^{l+p}},\ \tilde{d}_2 = \frac{2\left(\frac{m}{m+1}\right)^g - (1+\rho_w^2)}{{\ell_k}^g},\\
    & \Tilde{d}_3 = \frac{\left(\frac{m}{m+1}\right)^h-\frac{1+\rho_w^2}{2}-\frac{\alpha^2\gamma^2a_{33}}{(m+1)^{2p}}}{{\ell_k}^h}.
    \vspace{-2mm}
\end{split}
\end{equation*}
Note that $m^p \hspace{-.5mm}-\hspace{-.5mm}\frac{m^{p+l}}{(m+1)^l}$ decreases to 0 as $m\rightarrow\infty$ due to (\ref{eq:partiald1}) and L'Hopital's rule. There exists an $m$ satisfying the following conditions to ensure $\tilde{d}_1 >0$, $\tilde{d}_2 >0$, and $\tilde{d}_2 >0$ respectively: 
\begin{subequations} \label{eq:condm0p1}
    \begin{equation}
    \hspace{-4mm}\left\{\begin{array}{cl}
    \left(\frac{m}{m+1}\right)^l \hspace{-.5mm}>\hspace{-.5mm} 1-\alpha\gamma\mu & \text{ for } p=0, l>0\\
    m^p \hspace{-.5mm}-\hspace{-.5mm}\frac{m^{p+l}}{(m+1)^l} \hspace{-.5mm}<\hspace{-.5mm} \alpha\gamma\mu & \text{ for } 0\hspace{-.5mm}<\hspace{-.5mm}p\hspace{-.5mm}<\hspace{-.5mm}l\hspace{-.5mm}<\hspace{-.5mm} 1, m\hspace{-.5mm}\geq\hspace{-.5mm}\frac{p}{l-p}\\
    \alpha\gamma>\frac{l}{\mu} & \text{ for } 0<p\leq1\leq l
    \end{array}\right.,
    \end{equation}
    \vspace{-2mm}
    \begin{equation}
        \hspace{-5.5cm}\left(\frac{m}{m+1}\right)^g \hspace{-.5mm}>\hspace{-.5mm} \frac{1+\rho_w^2}{2},
    \end{equation}
    \begin{equation}
        \hspace{-4mm}\left\{\begin{array}{ll}
        \left(\hspace{-.5mm}\frac{m}{m+1}\hspace{-.5mm}\right)^h\hspace{-.5mm}>\hspace{-.5mm}\frac{1+\rho_w^2}{2}\hspace{-.5mm}+\hspace{-.5mm}\frac{\alpha^2\gamma^2a_{33}}{(m+1)^{2p}}, \alpha\gamma\hspace{-.5mm}\leq\hspace{-.5mm}\frac{1\hspace{-.5mm}-\hspace{-.5mm}\rho_w^2}{4\sqrt{2}Ld_I}    & \text{ for } p \hspace{-.5mm}=\hspace{-.5mm} 0\\
        \left(\hspace{-.5mm}\frac{m}{m+1}\hspace{-.5mm}\right)^h\hspace{-.5mm}>\hspace{-.5mm}\frac{1+\rho_w^2}{2}\hspace{-.5mm}+\hspace{-.5mm}\frac{\alpha^2\gamma^2a_{33}}{(m+1)^{2p}} & \text{ for } 0\hspace{-.5mm}<\hspace{-.5mm}p\hspace{-.5mm}\leq\hspace{-.5mm} 1
        \end{array}  \right..
    \end{equation}
\end{subequations}
Denote 
\vspace{-2mm}
\begin{equation} \label{eq:c}
\begin{split}
    & c_1\triangleq\alpha\gamma\mu \hspace{-.5mm}-\hspace{-.5mm} c(p,l), c_2\triangleq 2\left(\frac{m}{m \hspace{-.5mm}+\hspace{-.5mm} 1}\right)^g \hspace{-1.5mm}-\hspace{-.5mm} (1 \hspace{-.5mm}+\hspace{-.5mm} \rho_w^2),\\
    & c_3\triangleq \left(\frac{m}{m+1}\right)^h-\frac{1+\rho_w^2}{2}-\frac{\alpha^2\gamma^2a_{33}}{(m+1)^{2p}}.
\end{split}
\vspace{-2mm}
\end{equation}
Plugging $\Tilde{d}_1$, $\Tilde{d}_2$, and $\Tilde{d}_3$ into (\ref{eq:bdineq}) leads to 
\vspace{-2mm}
\begin{equation} \label{eq:bdineq0p}
    \begin{array}{ll}
        \hat{U} & \geq \frac{\alpha\gamma\tilde{a}_{12}}{c_1{\ell_k}^{g-l}}\hat{X} \hspace{-.5mm}+\hspace{-.5mm} \frac{b_1}{c_1{\ell_k}^{2q-p-l}}, 
        \hat{X} \geq \frac{a_{23}\alpha^2}{c_2{\ell_k}^{h-g}}\hat{Y} \hspace{-.5mm}+\hspace{-.5mm} \frac{b_2}{c_2{\ell_k}^{2q-g}},\\
        \hat{Y} & \geq \left(\frac{\alpha^2\gamma^4a_{31}}{c_3{\ell_k}^{2p+l-h}(\ell_k+1)^{2p}} + \frac{\gamma^2a_{34}c_{mp}^2}{c_3{\ell_k}^{2p+l-h}}\right)\hat{U}+\frac{\gamma^2b_{32}c_{mp}^2}{c_3{\ell_k}^{2p-h}}\\
        & \hspace{3.5mm}
        + \left( \frac{\gamma^2\tilde{a}_{32}}{c_3{\ell_k}^{g-h}(\ell_k+1)^{2p}} + \frac{\gamma^2a_{35}c_{mp}^2}{c_3{\ell_k^{2p+g-h}}}\right)\hat{X} + \frac{\tilde{b}_{31}}{c_3{\ell_k}^{2q-h}}.
        \vspace{-1mm}
    \end{array}
    \vspace{-2mm}
\end{equation}}\\
{The expressions on the {right-hand} side of the inequalities of (\ref{eq:bdineq0p}) 
are bounded by constants when the coefficients $l$, $h$, and $g$ satisfy the following linear inequalities:
\vspace{-2mm}
\begin{equation} \label{eq:lgh_0p1}
\begin{split}
    & 0< p < l \leq g \leq h \leq \min\{2q,2p\},\\
    & h - 2p \leq g \leq 2q, \ 
    h- 2p \leq l \leq 2q-p.
\end{split}
\end{equation}
or
\vspace{-2mm}
\begin{equation} \label{eq:lgh_p0}
\begin{split}
    & 0= p < l \leq g \leq h \leq 2q,\\
    & h - 2p \leq g \leq 2q, \ 
    h- 2p \leq l \leq 2q-p.
\end{split}
\end{equation}
When $q \geq p $, the solution set of $l$, $g$, and $h$ in (\ref{eq:lgh_0p1}) and (\ref{eq:lgh_p0}) is not empty. To get a tight error bound in (\ref{eq:UXYk}) and (\ref{eq:fgh_leq1}), we maximize the value of $l$, $g$, and $h$  and get the following solution:
\vspace{-2mm}
\begin{equation} \label{eq:lghvalue}
    \left\{\begin{array}{cl}
        l \hspace{-.5mm}=\hspace{-.5mm} \min\{2q\hspace{-.5mm}-\hspace{-.5mm}p,2p\},\ g \hspace{-.5mm}=\hspace{-.5mm} h =2\min\{q,p\} & \text{ if } 0 \hspace{-.5mm}<\hspace{-.5mm} p \hspace{-.5mm}\leq\hspace{-.5mm} 1 \\
        l \hspace{-.5mm}=\hspace{-.5mm} h \hspace{-.5mm}=\hspace{-.5mm} g \hspace{-.5mm}=\hspace{-.5mm} 2q & \text{ if } p \hspace{-.5mm}=\hspace{-.5mm} 0
    \end{array}\right..
\end{equation}}

Under (\ref{eq:bdineq0p}) and (\ref{eq:lghvalue}), 
we can set $e_0$, $e_1$, $e_2$, $e_3$, $e_4$, $e_5$, and $e_6$ in (\ref{eq:eineq}) as 
$e_0 = \frac{\alpha^2\gamma^4a_{31}}{c_3{m}^{4p+l-h}} + \frac{\gamma^2a_{34}c_{mp}^2}{c_3{m}^{2p+l-h}}$, 
$e_1 = \frac{\gamma^2\tilde{a}_{32}+\gamma^2a_{35}c_{mp}^2}{c_3{m^{2p}}}$, 
$e_2 = \frac{\gamma^2b_{32}c_{mp}^2}{c_3{m}^{2p-h}}+\frac{\tilde{b}_{31}}{c_3{m}^{2q-h}}$, 
$e_3=\frac{a_{23}\alpha^2}{c_2}$, 
$e_4=\frac{b_2}{c_2m^{2q-g}}$, $e_5=\frac{\alpha\gamma\tilde{a}_{12}}{c_1m^{g-l}}$, and $e_6=\frac{b_1}{c_1m^{2q-p-l}}$. 

Next, we give the condition guaranteeing (\ref{eq:condtodiscuss}), with
\begin{equation}
\begin{split}
    \frac{1}{e_3} \hspace{-.5mm}-e_0e_5 \hspace{-.5mm}-\hspace{-.5mm} e_1  \hspace{-.5mm}= & 
    \frac{c_2}{a_{23}\alpha^2} - \hspace{-1mm}
    \left(\hspace{-1mm} \frac{\gamma^2\tilde{a}_{32} \hspace{-.5mm}+\hspace{-.5mm} \gamma^2a_{35}c_{mp}^2}{c_3{m^{2p}}}\hspace{-1mm}\right)\\
    & \hspace{-.5mm}- \hspace{-1mm}\left(\hspace{-1mm}\frac{\alpha^3\gamma^5\tilde{a}_{12}a_{31}}{c_1c_3{m}^{4p}} {+ \frac{\alpha\gamma^3\tilde{a}_{12}a_{34}c_{mp}^2}{c_1c_3{m}^{2p}}}\hspace{-1mm}\right)\hspace{-1mm}.
\end{split}
\end{equation}
1) For $p = 0$ and $\alpha\gamma_0<\frac{2}{\mu+L}$, 
we have $c_{mp}=0$ and
$\frac{1}{e_3} \hspace{-.5mm}-e_0e_5 \hspace{-.5mm}-\hspace{-.5mm} e_1   = 1 \hspace{-.5mm}- \hspace{-.5mm}\left(\hspace{-1mm}\frac{\alpha^5\gamma^5\tilde{a}_{12}a_{23}a_{31}}{c_1c_2c_3} \hspace{-1mm}+\hspace{-1mm}\frac{a_{23}\alpha^2\gamma^2\tilde{a}_{32}}{c_2c_3}\hspace{-1mm}\right)\hspace{-1.5mm}>\hspace{-.5mm} c_4,$
with $c_4\hspace{-.5mm}\triangleq\hspace{-.5mm}1-{\alpha^5\gamma^5}\frac{a_{23}a_{31}L^2(3\mu+L)}{\mu n(\mu+L)c_1c_2c_3}-{\alpha^2\gamma^2\frac{32a_{23}d_I^2L^2{(\mu^2 + 2\mu L + 5L^2)}}{c_2c_3(1-\rho_w^2)(\mu+L)^2}}$.}
Note that $c_1$, $c_2$, $c_3$ monotonically increase with $m$, with their limits given by $\lim_{m\rightarrow\infty}c_1 = \alpha\gamma\mu$, $\lim_{m\rightarrow\infty}c_2 = 1-\rho_w^2$, and $\lim_{m\rightarrow\infty}c_3 = \frac{1-\rho_w^2}{2}$. Then as $m\rightarrow \infty$, $c_4$ monotonically increases to $\lim_{m\rightarrow\infty}c_4=1-c_{41}\alpha^4\gamma^4-c_{42}\alpha^2\gamma^2$ with $c_{41}\triangleq\frac{64d_I^2L^6({3\mu+L})(1+\rho_w^2)}{({\mu+L})(1-\rho_w^2)^4}$ and $c_{42}\triangleq\frac{64(1+\rho_w^2)d_I^2L^2(\mu^2+2\mu L+5L^2)}{({\mu+L}^2)(1-\rho_w^2)^4}$. Hence,
$\lim_{m\rightarrow\infty}c_4>0$ is guaranteed by 
\vspace{-2mm}
\begin{equation} \label{eq:alphagamma}
    \alpha\gamma<\sqrt{\frac{2}{c_{42}+\sqrt{c_{42}^2+4c_{41}}}}.
\end{equation}
Therefore, (\ref{eq:condtodiscuss}) is guaranteed by a large enough $m$ satisfying
\vspace{-2mm}
\begin{equation}  \label{eq:condm0p}
    c_4 > 0
\end{equation}
2) For $0<p\leq1$, noting that $c_1$, $c_2$, $c_3$ increase as $m$ increases, $c_{mp}$, $\tilde{a}_{12}$, and $\tilde{a}_{32}$ decrease as $m$ increases, and $\frac{1}{e_3} \hspace{-.5mm}-e_0e_5 \hspace{-.5mm}-\hspace{-.5mm} e_1$ increases to $\frac{(1-\rho_w^2)^2}{(1+\rho_w^2)\alpha^2}>0$ as $m\rightarrow\infty$, so there always exists an $m$ satisfying the following condition to ensure (\ref{eq:condtodiscuss}):
\vspace{-2mm}
\begin{equation} \label{eq:condm0p1_2}
\begin{split}
& \frac{\alpha^5\gamma^5\tilde{a}_{12}a_{23}\alpha^2a_{31}}{c_1c_2c_3m^{4p}}+\frac{a_{23}\alpha^2\gamma^2\tilde{a}_{32}}{c_2c_3m^{2p}} \\
& +\left(\frac{\gamma^3\alpha^3\tilde{a}_{12}a_{23}a_{34}}{c_1c_2c_3m^{2p}} + \frac{\gamma^2\alpha^2a_{23}{a}_{35}}{c_2c_3m^{2p}}\right)c_{mp}^2< 1.
\end{split}
\end{equation}

{
In summary, combining (\ref{eq:alphagamma0}), (\ref{eq:condm0p1}), (\ref{eq:lghvalue}),  (\ref{eq:alphagamma}), and (\ref{eq:condm0p}) gives the conditions for $\alpha$, $\gamma$ and $m$ for the case $p=0$, and combining (\ref{eq:alphagamma0}), (\ref{eq:condm0p1}), (\ref{eq:lghvalue}), and (\ref{eq:condm0p1_2}) gives the conditions for $\alpha$, $\gamma$ and $m$ for the case $0<p\leq1$. Combining (\ref{eq:errorevo}), (\ref{eq:fgh_leq1}), and (\ref{eq:lghvalue}) gives the convergence properties for the cases $p=0$ and $0<p\leq1$.
}

\textbf{Case II:} For $p>1$, we choose $f_k$, $g_k$, and $h_k$ in the following form: 
\vspace{-4mm}
\begin{equation} \label{eq:fgh_leq2}
    f_k = 1+\frac{1}{(m+k)^p}, g_k = \frac{1}{(m+k)^g}, h_k = \frac{1}{(m+k)^h},
\end{equation}
where $g$ and $h$ are some positive variables depending on $p$ and $q$. Plugging (\ref{eq:fgh_leq2}) into (\ref{eq:d}) yields
\vspace{-4mm}
\begin{equation}
\begin{split}
    &d_1 
    = \frac{\alpha\gamma\mu}{{\ell_k}^p} + \frac{\left(1-\frac{{\ell_k}}{{\ell_{k}+1}}\right)^p-1}{{\ell_k}^p} + \frac{\alpha\gamma\mu}{{\ell_k}^{2p}}> \frac{\alpha\gamma\mu}{{\ell_k}^p} - \frac{ p}{{\ell_k}^{p+1}}
,\\
\vspace{-1mm}
    &d_2
    = \frac{2\left(\frac{\ell_{k}}{{\ell_{k}+1}}\right)^g - (1+\rho_w^2)}{{\ell_k}^g},\\
    \vspace{-1mm}
    & d_3= \frac{\left(\frac{\ell_{k}}{{\ell_{k}+1}}\right)^h-(\frac{1+\rho_w^2}{2}+\frac{\alpha^2\gamma^2a_{33}}{(\ell_{k}+1)^{2p}})}{{\ell_k}^h}.
\end{split}
    \vspace{-10mm}
\end{equation}
Select $\tilde{d}_1$, $\tilde{d}_2$, and $\tilde{d}_3$ as 
\vspace{-6mm}
\begin{equation} 
    \begin{split}
        &
        \tilde{d}_1 = \frac{\alpha\gamma\mu-\frac{p}{m}}{{\ell_k}^p}, 
        \tilde{d}_2 = \frac{2\left(\frac{m}{m+1}\right)^g \hspace{-1mm}-\hspace{-.5mm} (1+\rho_w^2)}{{\ell_k}^g},\\
        & \tilde{d}_3 = \frac{\left(\frac{m}{m+1}\right)^h\hspace{-1mm}-\hspace{-.5mm}\frac{1+\rho_w^2}{2}\hspace{-.5mm}-\hspace{-.5mm}\frac{\alpha^2\gamma^2a_{33}}{(m+1)^{2p}}}{{\ell_k}^h}.
        \vspace{-2mm}
    \end{split}
\end{equation}
$\tilde{d}_1>0$ holds when
\vspace{-4mm}
\begin{equation} \label{eq:condmp0}
    m>\frac{p}{\alpha\gamma\mu}
    \vspace{-3mm}
\end{equation}
holds. $\tilde{d}_2>0$ and $\tilde{d}_3>0$ hold when $m$ is large enough to satisfy
 \vspace{-3mm}
 \begin{equation} \label{eq:condmp0_2}
     \left(\frac{m}{m+1}\right)^g > \frac{1+\rho_w^2}{2} \hspace{-.5mm}+\hspace{-.5mm}\frac{\alpha^2\gamma^2a_{33}}{(m + 1)^{2p}}.
     \vspace{-3mm}
 \end{equation}
Denote
\vspace{-3mm}
\begin{equation}
\begin{split}
    & c_1\triangleq\alpha\gamma\mu-\frac{p}{m}, c_2\triangleq 2\left(\frac{m}{m+1}\right)^g - (1+\rho_w^2),\\
    & c_3\triangleq \left(\frac{m}{m+1}\right)^h-\frac{1+\rho_w^2}{2}-\frac{\alpha^2\gamma^2a_{33}}{(m+1)^{2p}},
\end{split}
\vspace{-4mm}
\end{equation}
then (\ref{eq:bdineq}) gives 
\vspace{-4mm}
\begin{equation} \label{eq:bdineqp0}
    \begin{array}{ll}
        \hat{U} & \geq \frac{\alpha\gamma\tilde{a}_{12}}{c_1{\ell_k}^{g}}\hat{X} + \frac{b_1}{c_1{\ell_k}^{2q-p}}, 
        \hat{X} \geq \frac{a_{23}\alpha^2}{c_2{\ell_k}^{h-g}}\hat{Y}+\frac{b_2}{c_2{\ell_k}^{2q-g}},\\
        \hat{Y} & \geq \left(\frac{\alpha^2\gamma^4a_{31} }{c_3{\ell_k}^{2p-h}(\ell_k+1)^{2p}}+ \frac{\gamma^2a_{34}c_{mp}^2}{c_3{\ell_k}^{2p-h}} \right)\hspace{-1mm}\left(1+\frac{1}{{\ell_k}^p}\right)\hat{U}\\
        & \hspace{3.5mm} + \left(\frac{\gamma^2\tilde{a}_{32}}{c_3{\ell_k}^{g-h}(\ell_k+1)^{2p}} +\frac{\gamma^2{a}_{35}c_{mp}^2}{c_3{\ell_k}^{2p+g-h}} \right)\hat{X}\\
        & \hspace{3.5mm} + \frac{\tilde{b}_{31}}{c_3{\ell_k}^{2q-h}}+\frac{\gamma^2b_{32}c_{mp}^2}{c_3{\ell_k}^{2p-h}}.
    \end{array}
\end{equation}

Under $q\geq\frac{p}{2}$, the expressions on the right-hand side of the inequalities of (\ref{eq:bdineqp0}) are bounded by constants when the coefficients $l$, $h$, and $g$ in (\ref{eq:fgh_leq2}) satisfy the following inequalities:
\vspace{-2mm}
\begin{equation} 
\begin{split}
    0 \leq g \leq h \leq \min\{2q,2p\},\ 
    h - 2p \leq g \leq 2q.\\
\end{split}
\end{equation}
We maximize the values of $g$ and $h$ and get the following solution:
\vspace{-2mm}
\begin{equation} \label{eq:lghvalue2}
    g = h =\min\{2q,2p\}.
\end{equation}

Under (\ref{eq:bdineqp0}) and (\ref{eq:lghvalue2}), we can set $e_0$, $e_1$, $e_2$, $e_3$, $e_4$, $e_5$, and $e_6$ in (\ref{eq:eineq}) as $e_0=\left(\frac{\alpha^2\gamma^4a_{31} }{c_3{m}^{4p-h}}+ \frac{\gamma^2a_{34}c_{mp}^2}{c_3{m}^{2p-h}} \right)\hspace{-1mm}\left(1+\frac{1}{{m}^p}\right)$, $e_1 = \frac{\gamma^2\tilde{a}_{32}}{c_3m^{2p}} {+\frac{\gamma^2{a}_{35}c_{mp}^2}{c_3{m}^{2p}}}$, $e_2 = \frac{\tilde{b}_{31}}{c_3{m}^{2q-h}}+\frac{\gamma^2b_{32}c_{mp}^2}{c_3{m}^{2p-h}}$, $e_3=\frac{a_{23}\alpha^2}{c_2}$, $e_4=\frac{b_2}{c_2m^{2q-g}}$, $e_5=\frac{\alpha\gamma\tilde{a}_{12}}{c_1m^{g}}$, and $e_6=\frac{b_1}{c_1m^{2q-p}}$.

Next, we give the condition guaranteeing (\ref{eq:condtodiscuss}), with
\begin{equation}
\begin{split}
    \frac{1}{e_3} \hspace{-.5mm}-\hspace{-.5mm} e_0e_5 \hspace{-.5mm}-\hspace{-.5mm} e_1  \hspace{-1mm} = & 
    \frac{c_2}{a_{23}\alpha^2} \hspace{-.5mm}-\hspace{-1mm} 
    \left(\hspace{-1mm}\frac{\gamma^2\tilde{a}_{32}}{c_3m^{2p}} {+\frac{\gamma^2{a}_{35}c_{mp}^2}{c_3{m}^{2p}}} \hspace{-1mm}\right)\hspace{-1mm}-\hspace{-1.5mm}\left(\hspace{-1mm}1 \hspace{-1mm}+\hspace{-1mm}\frac{1}{{m}^p}\hspace{-1mm}\right)\\
    & 
    \hspace{-.5mm}\times\hspace{-1mm}\left(\hspace{-1mm}\frac{\alpha^3\gamma^5\tilde{a}_{12}a_{31} }{c_1c_3{m}^{4p}} {+ \frac{\alpha\gamma^3\tilde{a}_{12}a_{34}c_{mp}^2}{c_1c_3{m}^{2p}} }\hspace{-1mm}\right)\hspace{-1mm}.
\end{split}
\end{equation}
For $p>1$, note that $c_1$, $c_2$, $c_3$ increase as $m$ increases, and $c_{mp}$, $\tilde{a}_{12}$, and $\tilde{a}_{32}$ decrease as $m$ increases, and $\frac{1}{e_3} \hspace{-.5mm}-e_0e_5 \hspace{-.5mm}-\hspace{-.5mm} e_1$ increases to $\frac{(1-\rho_w^2)^2}{(1+\rho_w^2)\alpha^2}>0$ as $m\rightarrow\infty$, so (\ref{eq:condtodiscuss}) is guaranteed by a large enough $m$ satisfying
\vspace{-2mm}
\begin{equation} \label{eq:condmp0_3}
    \begin{split}
    & \frac{\alpha^5\gamma^5\tilde{a}_{12}a_{23}\alpha^2a_{31}}{c_1c_2c_3m^{4p}}\hspace{-1mm}\left(\hspace{-.5mm}1\hspace{-.5mm}+\hspace{-.5mm}\frac{1}{m^p}\hspace{-.5mm}\right)\hspace{-.5mm}+\hspace{-.5mm}\frac{a_{23}\alpha^2\gamma^2\tilde{a}_{32}}{c_2c_3m^{2p}} \hspace{-.5mm}\\
    &
    + \hspace{-1mm} \left(  \hspace{-1mm} \frac{\alpha^3\gamma^3\tilde{a}_{12}a_{23}a_{34} }{c_1c_2c_3{m}^{2p}} \hspace{-1mm}
    \left(\hspace{-.5mm} 1 \hspace{-.5mm}+\hspace{-.5mm}\frac{1}{{m}^p} \hspace{-.5mm}\right)    \hspace{-1mm}+\hspace{-.5mm} \frac{\alpha^2\gamma^2a_{23}{a}_{35}}{c_2c_3{m}^{2p}}\hspace{-1mm}\right)\hspace{-1mm}c_{mp}^2 \hspace{-1mm}
    <\hspace{-.5mm} 1.
    \end{split}
\end{equation}
In summary, combining (\ref{eq:alphagamma0}), (\ref{eq:condmp0}), (\ref{eq:condmp0_2}), (\ref{eq:lghvalue2}), and (\ref{eq:condmp0_3}) gives the conditions for $\alpha$, $\gamma$ and $m$ for the case $p>1$. Combining (\ref{eq:errorevo}), (\ref{eq:fgh_leq2}), and (\ref{eq:lghvalue2}) gives the convergence properties for the case $p>1$.

%% file: 64_appendix4.tex
\section{Proof of Theorem~\ref{thmrho}} \label{sec:rhow}

{Using the expressions of $\theta_1$ and $\theta_2$ in Theorem~\ref{thmConv0p}, we have the optimization error $\theta\triangleq 2n\theta_1 + 2\theta_2$ as
\vspace{-2mm}
\begin{equation*}
    \theta =
    \frac{e_1 T_w^4 \hspace{-.5mm}+\hspace{-.5mm} e_2 T_w^3 \hspace{-.5mm}+\hspace{-.5mm} e_3 T_w^2 \hspace{-.5mm}+\hspace{-.5mm} e_4 T_w \hspace{-.5mm}+\hspace{-.5mm} e_5}{c_1 T_w^4 \hspace{-.5mm}+\hspace{-.5mm} c_2 T_w^2 \hspace{-.5mm}+\hspace{-.5mm} c_3 T_w \hspace{-.5mm}+\hspace{-.5mm} c_4}\rho({W_o})^2
\end{equation*}
with 
$e_1 = \frac{n\alpha^2}{4}\sigma_\eta^2 + \frac{n}{4}\sigma_\xi^2$,
$e_2
= \left( \frac{n\alpha{d_I^2}L^2}{\mu} + \frac{\alpha\mu}{2}n{d_I^2}\right)\sigma_\xi^2$,
$e_3 = \left(-8n{d_I^2}\alpha^4L^2\right)\sigma_\eta^2 + \left(-8n{d_I^2}\alpha^2L^2\right)\sigma_\xi^2$,
$e_4 = \left(-\frac{32n\alpha^3 {d_I^4}L^4}{\mu}-16n\alpha^3\mu{d_I^4}L^2\right)\sigma_\xi^2$, and
$e_5= \left(\frac{96n\alpha^3 {d_I^2}L^2}{\mu} -\right.$
$\left. 128n \alpha^4 {d_I^2}L^2 + 64n\alpha^6{d_I^2}L^4 + 48\mu n \alpha^3{d_I^2}\right)\sigma_\eta^2 + \left(\frac{16n\alpha^3{d_I^2}L^4}{\mu} \right.$
$\left. -128n \alpha^2 {d_I^2}L^2 + {64n\alpha^4 {d_I^2}L^4} + 8\mu n\alpha^3{d_I^2}L^2\right)\sigma_\xi^2$. The parameters $T_w$ and $c_1 - c_4$ are given in Theorem~\ref{thmConv0p}.}

According to the definition of $T_w$, we can prove $\frac{\partial\theta}{\partial{\rho_w}}>0$ by proving $\frac{\partial\theta}{\partial{T_w}}<0$ where
{
\vspace{-2mm}
\begin{equation*}
    \frac{\partial\theta}{\partial T_w} \hspace{-.5mm}=\hspace{-.5mm} 
    \frac{-\hspace{-.5mm}f_1 T_w^6 \hspace{-.5mm}+\hspace{-.5mm} f_3 T_w^4 \hspace{-.5mm} - \hspace{-.5mm} f_4 T_w^3 \hspace{-.5mm}-\hspace{-.5mm} f_5 T_w^2 \hspace{-.5mm}+\hspace{-.5mm} f_6 T_w \hspace{-.5mm}+\hspace{-.5mm} f_7} {(c_1 T_w^4 \hspace{-.5mm}+\hspace{-.5mm} c_2 T_w^2 \hspace{-.5mm}+\hspace{-.5mm} c_3)^2}\hspace{-.5mm}\rho({W_o})^2\hspace{-.5mm}
\end{equation*}
with $f_1$, $f_3$, $f_4$, $f_5$, $f_6$, and $f_7$ being some positive constants independent of $T_w$:
$f_1 = \left( \frac{n{d_I^2}L^2}{4} + \frac{\mu^2n{d_I^2}}{8}\right)\sigma_\xi^2{\alpha^2}$, 
$f_3 =\left( {16n {d_I^4}L^4}+8 n\mu^2{d_I^4}L^2\right)\sigma_\xi^2\alpha^4$,
$f_4 = \left(\frac{128n\alpha^7L^6}{\mu} \right.$
$\left.+ 64n\mu\alpha^7L^4 +{96n\alpha^4 L^2}  + 48 n\mu^2 \alpha^4\right){d_I^2}\sigma_\eta^2
+ \left(\frac{128n\alpha^5 L^6}{\mu} \right.$
$\left.+ {64n\mu\alpha^5 L^4}+ {16n\alpha^4L^4} + 8\mu^2 n\alpha^4L^2 \right){d_I^2}\sigma_\xi^2$,
$f_5 = \left( \frac{384nL^8\alpha^6}{\mu^2}  + {192n\alpha^6 L^6}+ {256n {d_I^2}L^6}\alpha^6 + 128n\mu^2{d_I^2}L^4\alpha^6\right.$
$\left. + {384nL^4}\alpha^4 + 192n\alpha^4 \mu^2 L^2\right){d_I^4}\sigma_\xi^2$,
$f_6 = 16\left( \frac{128n\alpha^9L^8}{\mu}\right.$
$\left. + 64n\mu\alpha^9L^6  + 96n\alpha^6L^4+ 48\mu^2 n \alpha^6L^2 \right){d_I^4}\sigma_\eta^2
+ \left( \frac{128n\alpha^7L^8}{\mu}\right.$ 
$\left.\hspace{-.5mm}+64n\mu\alpha^7\hspace{-.5mm}L^6 \hspace{-.5mm}+\hspace{-.5mm} 16n\alpha^6\hspace{-.5mm}L^6\hspace{-.5mm}+\hspace{-.5mm} 8\mu^2 n\alpha^6\hspace{-.5mm}L^4 \right)\hspace{-.5mm}16{d_I^4}\sigma_\xi^2$,
$f_7 \hspace{-.5mm}=\hspace{-.5mm} 16{d_I^6}\sigma_\xi^2$
\\
$\times\hspace{-.5mm}\left(\frac{256n\alpha^8L^{10}}{\mu^2}+{128n\alpha^8 L^8} + {256n\alpha^6L^6}+ 128n\alpha^6 \mu^2L^4\right)$.}
{When
\vspace{-2mm}
\begin{equation} \label{eq:alpha_thmrhow}
    \alpha^2\leq \frac{T_w^2}{64d_I^2L^2},
\end{equation}we have $f_7< f_7\frac{3T_w^2}{32d_I^2L^2\alpha^2}$, $f_6T_w< f_6\frac{T_w^2}{16d_I^2L^2\alpha^2}T_w=f_4T_w^3$, and $f_3T_w^4\leq f_3\frac{T_w^2}{64d_I^2L^2\alpha^2}T_w^4=f_1T_w^6$. Hence, we have $\frac{\partial\theta}{\partial T_w}<0$ holds as 
\vspace{-2mm}
\begin{equation*}
\begin{split}
    \frac{\partial\theta}{\partial T_w} & <
    \frac{ \left(- f_5 + f_7 \frac{d_I^2L^2}{16\alpha^2}\right)T_w^2\rho({W_o})^2}{(c_1 T_w^4 \hspace{-.5mm}+\hspace{-.5mm} c_2 T_w^2 \hspace{-.5mm}+\hspace{-.5mm} c_3)^2}\\
    & = \frac{\left(  - {256nL^6} - 128n\mu^2L^4\right)\alpha^6{d_I^6}\sigma_\xi^2T_w^2\rho({W_o})^2}{(c_1 T_w^4 \hspace{-.5mm}+\hspace{-.5mm} c_2 T_w^2 +\hspace{-.5mm} c_3)^2}<0.
    \vspace{-2mm}
\end{split}
\end{equation*}
Combining the conditions of step size $\alpha$ in Theorem~\ref{thmConv0p} and (\ref{eq:alpha_thmrhow}) gives the step size condition in the theorem.
Further using the relationship $T_w=1-\rho_w^2$, we have $\frac{\partial \theta}{\partial \rho_w} > 0$.}

From the expressions of $\theta_1$ and $\theta_2$ in (\ref{eq:theta1}), it can be seen that both $\theta_1$ and $\theta_2$ increase with an increase in $\rho(W_o)$. Therefore, $\theta$ also increases with an increase in $\rho(W_o)$.